\documentclass[a4paper,11pt]{article}
\title{\bf{Moduli stacks and invariants of semistable objects on 
K3 surfaces}
}
\date{}
\author{Yukinobu Toda}

\usepackage{makeidx}

\usepackage{latexsym}
\usepackage{amscd}
\usepackage{amsmath}
\usepackage{amssymb}

\usepackage{amsthm}
\usepackage{float}
\usepackage[all]{xy}
\usepackage{graphicx}

\usepackage{macros}

\usepackage{array}
\usepackage{amscd}
\usepackage[all]{xy}
\usepackage{makeidx}
\usepackage{latexsym}
\DeclareFontFamily{U}{rsfs}{%
\skewchar\font127}
\DeclareFontShape{U}{rsfs}{m}{n}{%
<-6>rsfs5<6-8.5>rsfs7<8.5->rsfs10}{}
\DeclareSymbolFont{rsfs}{U}{rsfs}{m}{n}
\DeclareSymbolFontAlphabet
{\mathrsfs}{rsfs}
\DeclareRobustCommand*\rsfs{%
\@fontswitch\relax\mathrsfs}
\theoremstyle{plain}
\newtheorem{thm}{Theorem}[section]
\newtheorem{prop}[thm]{Proposition}
\newtheorem{lem}[thm]{Lemma}
\newtheorem{defi}[thm]{Definition}
\newtheorem{rmk}[thm]{Remark}
\newtheorem{cor}[thm]{Corollary}
\newtheorem{prob}[thm]{Problem}
\newtheorem{step}{Step}
\newtheorem{sstep}{Step}
\newtheorem{ssstep}{Step}

\newtheorem{claim}[thm]{Claim}

\newtheorem{prop-defi}[thm]{Proposition-Definition}
\newtheorem{problem}[thm]{Problem}
\newtheorem{assum}[thm]{Assumption}
\newtheorem{conj}[thm]{Conjecture}
\newtheorem{exam}[thm]{Example}
\begin{document}
\maketitle

\begin{abstract}
For a K3 surface $X$ and its bounded derived category of 
coherent sheaves $D(X)$, we have the notion of stability 
conditions on $D(X)$ in the sense of T.~Bridgeland. 
In this paper, we show that the moduli stack of semistable 
objects in $D(X)$ with a fixed numerical class and a
phase is represented by an Artin stack of finite type over 
 $\mathbb{C}$. Then following 
D.~Joyce's work, we introduce the invariants counting semistable 
objects in $D(X)$, and show that the invariants are independent
of a choice of a stability condition. 
\end{abstract}

\section{Introduction}
The work of this paper is motivated by D.~Joyce's recent 
works~\cite{Joy1}, \cite{Joy2}, \cite{Joy3}, \cite{Joy4}, \cite{Joy},
especially~\cite[Conjecture 6.25]{Joy4}
on the counting invariants of semistable objects on K3 surfaces 
or abelian surfaces. Such invariants are expected to 
produce automorphic functions on the space of stability 
conditions in the sense of T.~Bridgeland~\cite{Brs1}. 

\subsection{Stability conditions}
Let $X$ be a smooth projective variety over $\mathbb{C}$, 
$\Coh(X)$ the abelian category of coherent sheaves on $X$, 
and $D(X)$ the bounded derived category of $\Coh (X)$. 
For an ample divisor $\omega$ on $X$, there is a notion of 
$\omega$-Gieseker stability 
on $\Coh (X)$, and the 
moduli spaces of semistable sheaves have been studied in detail up to 
now~\cite{Hu}. 
The notion of stability conditions on a triangulated category $\tT$ 
(especially including the case of $\tT=D(X)$) was introduced by T.~Bridgeland~\cite{Brs1} motivated by M.Douglas's
 $\Pi$-stability~\cite{Dou1}, \cite{Dou2}.
Roughly it consists of data $\sigma =(Z, \pP)$, 
$$Z\colon K(\tT) \lr \mathbb{C},\quad \pP (\phi)\subset \tT,$$
where $Z$ is a group homomorphism and $\pP(\phi)$ is a full subcategory 
for each $\phi \in \mathbb{R}$, and these data satisfy some axiom. 
(See Definition~\ref{general} below.)
Then Bridgeland~\cite{Brs1} showed that the 
set of good stability conditions has a structure 
of a complex manifold.
When $\tT=D(X)$, the space of stability conditions 
$\Stab (X)$ carries a map,
$$\zZ \colon \Stab (X) \lr \nN (X)_{\mathbb{C}}^{\ast},$$
where $\nN (X)=K(X)/\equiv$ is a \textit{numerical Grothendieck group}.
(See Definition~\ref{numGro}.) 
The precise descriptions of the space $\Stab (X)$ have been 
studied in the articles~\cite{Brs2}, \cite{Brs3}, 
\cite{Brs4}, \cite{Tho}, \cite{Oka}, \cite{Mac}, \cite{IUU}, \cite{Ber},
\cite{Tst}, \cite{Tst2}. In particular when $X$ is a K3 surface or 
an abelian surface, Bridgeland~\cite{Brs2} described 
$\Stab ^{\ast}(X)$, one 
of the connected components of $\Stab (X)$, as a covering space 
over a certain open subset $\pP _0 ^{+}(X)\subset \nN (X)^{\ast}_{\mathbb{C}}$, and related its Galois group to the group of autoequivalences of 
$D(X)$. 

In general
when $X$ is
a Calabi-Yau manifold, 
it is expected that the space $\Stab (X)$
describes the so called ``stringy K$\ddot{\textrm{a}}$hler moduli space".
More precisely Bridgeland
 conjectures in~\cite{Brs6} that the double quotient space, 
\begin{align}\label{autfu}
\Auteq D(X) \backslash \Stab (X)/ \mathbb{C},\end{align}
contains the stringy K$\ddot{\textrm{a}}$hler moduli space 
$\mM _K (X)$, 
which is in a mirror side $\hat{X}$, isomorphic to the moduli space of the 
complex structures $\mM _{C}(\hat{X})$.
When
$X$ is an elliptic curve, $\mM _C(\hat{X})$ is 
nothing but the modular curve, and
we have the following complete picture~\cite{Brs1}, 
$$\Auteq D(X) \backslash \Stab (X)/ \mathbb{C}\cong \mM _C (\hat{X})
=\hH / \SL (2, \mathbb{Z}),$$
where $\hH\subset \mathbb{C}$ is the upper half plane. 
On the space $\hH$, several  
automorphic functions (Eisenstein series, $j$-invariant)
have been studied. 
Thus it is interesting to construct
automorphic functions on the 
space $\Stab (X)$, purely from the categorical data of 
$D(X)$, and compare the classical theory in the 
mirror side.   

\subsection{Counting invariants of semistable sheaves}
D.~Joyce's recent
works~\cite{Joy1}, \cite{Joy2}, \cite{Joy3}, \cite{Joy4}, \cite{Joy}
are attempts to introduce some structures on the space 
$\Stab (X)$, such as Frobenius structures or automorphic functions. 
However for several technical reasons, his arguments work
only on the space of stability conditions on an abelian category. 
What we are interested in this paper 
is the work~\cite{Joy4}, where D.~Joyce studies certain 
counting invariants
of semistable sheaves on a K3 surface $X$. 
We denote $C(X)\subset \nN(X)$ the image of $\Coh (X) \to \nN (X)$, 
and let
$\alpha \in C(X)$ be a numerical class
 and $\Lambda$ a $\mathbb{Q}$-algebra.
We consider a motivic invariant,
\begin{align}\label{motivic}
\Upsilon \colon (\mbox{quasi-projective varieties}) \lr
\Lambda.\end{align}
As an example, one can take $\Lambda =\mathbb{Q}(z)$ and 
$\Upsilon(Y)$ to be the virtual Poincare polynomial of $Y$.  
Using $\Upsilon$, 
D.~Joyce~\cite{Joy4} constructs an invariant
  $\hat{I}^{\alpha}(\omega)\in \Lambda$
  which counts $\omega$-Gieseker semistable sheaves of 
  numerical type $\alpha$, and its 
  weighted counting
 \begin{align}\label{weigh}\hat{J}^{\alpha}(\omega)=
\sum _{\alpha _1 + \cdots +\alpha _n=\alpha}
l^{-\sum _{j>i}\chi (\alpha _j, \alpha _i)}
\frac{(-1)^{n-1}(l-1)}{n}\prod _{i=1}^n \hat{I}^{\alpha _i}
(\omega) \in \Lambda.
\end{align}
Here $\alpha _i \in C(X)$ has the same reduced Hilbert polynomial 
with $\alpha$ and $l=\Upsilon (\mathbb{A}^1)\in \Lambda$.
Then Joyce~\cite{Joy4} showed that $\hat{J}^{\alpha}(\omega)$
does not depend on a choice of $\omega$, so one can denote it by
$\hat{J}^{\alpha}\in \Lambda$. 

The purpose of this paper is to translate the above work
into the context of Bridgeland's stability conditions.
As we see below, it is related to the automorphic functions on 
the space of stability conditions. 
Based on the results in~\cite{Joy4},
 D.~Joyce proposes the following 
conjecture. 
\begin{conj}\label{kotae}\emph{\bf{\cite[Conjecture 6.25]{Joy4}}} 
Let $X$ be a K3 surface or an abelian surface. 
For $\sigma \in \Stab ^{\ast}(X)$ and 
$\alpha \in \nN (X)$, there is  
$J^{\alpha}(\sigma)\in \Lambda$, a certain weighted counting of 
$\sigma$-semistable objects of numerical type $\alpha$,
such that 

(i) $J^{\alpha}(\sigma)$ does not depend on a choice of $\sigma$. 
Hence we can write it $J^{\alpha}\in \Lambda$.

(ii) If $\alpha \in C(X)$, then $J^{\alpha}=\hat{J}^{\alpha}$. 
\end{conj} 
Suppose for instance Conjecture~\ref{kotae} is true. 
 Let $\Auteq ^{\ast}D(X)$ be the group of autoequivalences on $D(X)$
which preserve the component $\Stab ^{\ast}(X)$. 
Then the property (i) of Conjecture~\ref{kotae} implies 
that 
$J^{\alpha}=J^{\Phi _{\ast}
\alpha}$ for $\Phi \in \Auteq ^{\ast}D(X)$.  
(See Corollary~\ref{yume} below.)
Based on this observation, Joyce~\cite{Joy4} suggests that the 
map (ignoring convergence)
\begin{align}\label{conver}
\Stab ^{\ast}(X) \ni \sigma =(Z, \pP) \longmapsto \sum _{\alpha \in \nN (X)\setminus \{0\}} \frac{J^{\alpha}}{Z(\alpha)^k}
\in \Lambda \otimes _{\mathbb{Q}}\mathbb{C}, \end{align}
for $k\in \mathbb{Z}$
would give a holomorphic function on $\Stab ^{\ast}(X)$ which is 
invariant under the action of $\Auteq ^{\ast}D(X)$,
i.e. automorphic function on $\Stab ^{\ast}(X)$. 
Our goal is the following. 
\begin{thm}\label{mainth}
Conjecture~\ref{kotae} is true.
\end{thm}
As stated in~\cite{ICM}, \cite{Joy4},
it is interesting to compare the formula (\ref{conver}) with 
the work of Borcherds~\cite{Borch} on the product expansions of the 
automorphic forms. 
\subsection{Moduli problems}
The first issue in attacking Conjecture~\ref{kotae}
is to develop the moduli theory of semistable 
objects in the sense of Bridgeland. 
The moduli theory of objects in $D(X)$
is studied in 
some articles~\cite{Inaba}, \cite{Inmo}, \cite{LIE}, \cite{AP}. 
In the recent work of Inaba~\cite{Inmo}, he constructs 
some nice moduli spaces of complexes, using the notion of 
ample sequences.  
However the relationship between Bridgeland's stability 
conditions~\cite{Brs1} and Inaba's stability conditions using 
ample sequences~\cite{Inmo} is not clear. 
On the other hand,  for our purpose we do not require the 
 moduli spaces to have good properties, (projective, fine, etc).
In fact 
we only need it to be an Artin stack of finite type. 
Thus in Section~\ref{stack}, we work over $D(X)$ for an 
arbitrary smooth projective variety $X$ and 
establish the general arguments to guarantee the 
moduli stacks to be Artin stacks of finite type. 

In Section~\ref{stack}, the work of Lieblich~\cite{LIE} 
would help us. Let $\mM$ be the moduli stack of 
objects $E\in D(X)$ which satisfies $\Ext ^{<0}(E,E)=0$. 
Then he showed that $\mM$ is an Artin stack of locally 
finite type over $\mathbb{C}$. 
For $\alpha \in \nN (X)$, $\phi \in \mathbb{R}$ and $\sigma =(Z, \pP) \in 
\Stab (X)$, we study the substack,
$$\mM ^{(\alpha, \phi)}(\sigma) \subset \mM,$$
which is the moduli stack of $E\in \pP (\phi)$ and of numerical type $\alpha$.
At least we have to resolve the following two problems,
addressed by~\cite{AP}. 
\begin{itemize}
\item {\bf Generic flatness problem}

Let $\aA =\pP ((0,1])\subset D(X)$ and $\eE \in D(X\times S)$
a family of objects in $D(X)$. Then is the locus $s\in S$
on which $\eE _s \in \aA$ an open subset of $S$?

\item {\bf Boundedness problem}

Does the set of objects $E\in \pP (\phi)$ of numerical type $\alpha$ form a 
bounded family?
\end{itemize} 
Let $\Stab ^{\ast}(X)$ be one of the connected components of the 
space $\Stab (X)$, and suppose it satisfies
the Assumption~\ref{nbd} below. 
Especially we require a certain subset $\vV \subset \Stab ^{\ast}(X)$,
which in several examples
can be taken to be the so called 
\textit{neighborhoods of the large volume limits}. 
The main theorem in Section~\ref{stack} is the following. 
\begin{thm}\label{main2}\emph{\bf{[Theorem~\ref{condition}]}}
Assume the generic flatness problem and the boundedness problem 
are true for any $\sigma \in \vV$.
Then for any  
$\sigma \in \Stab ^{\ast}(X)$, $\alpha \in \nN (X)$ and 
$\phi \in \mathbb{R}$, the stack 
$\mM ^{(\alpha, \phi)}(\sigma)$ is an Artin stack of finite 
type over $\mathbb{C}$. 
\end{thm}

In Section~\ref{K3}, we check that the assumption in Theorem~\ref{main2}
is satisfied when $X$ is a K3 surface or an abelian surface, 
and $\Stab ^{\ast}(X)$ is the connected component 
described in~\cite{Brs3}.
Thus we obtain the following. 
\begin{thm}\label{main3}\emph{\bf{[Theorem~\ref{sekai}]}}
Let $X$ be a K3 surface or an abelian surface. Then 
for any 
$\sigma \in \Stab ^{\ast}(X)$, 
$\alpha \in \nN (X)$ and $\phi \in \mathbb{R}$,
 the stack 
$\mM ^{(\alpha, \phi)}(\sigma)$ is an Artin stack of 
finite type over $\mathbb{C}$. 
\end{thm}

\subsection{Counting invariants of semistable objects}
The next step is to study the invariant
determined by the moduli 
stack $\mM ^{(\alpha, \phi)}(\sigma)$.
Given data (\ref{motivic}),   
we introduce $J^{\alpha}(\sigma)\in \Lambda$
for $\alpha \in \nN (X)$ and 
$\sigma \in \Stab ^{\ast}(X)$ in 
a completely similar way of $\hat{J}^{\alpha}(\omega)$. 
Then translating the arguments in~\cite{Joy4} to the 
context of Bridgeland's stability conditions, we show the following
in Section~\ref{counting}.  
\begin{thm}\label{main3}\emph{\bf{[Theorem~\ref{inde}]}}
The invariant $J^{\alpha}(\sigma)\in \Lambda$ does not 
depend on a choice of $\sigma \in \Stab ^{\ast}(X)$. 
\end{thm}

Thus we may write $J^{\alpha}(\sigma)=J^{\alpha}$.
Finally in Section~\ref{compare}, we compare
$J^{\alpha}$ and $\hat{J}^{\alpha}$. 
\begin{thm}\emph{\bf{[Theorem~\ref{goal}]}}\label{main4}
For $\alpha \in C(X)$, we have $J^{\alpha}=\hat{J}^{\alpha}$. 
\end{thm} 
By Theorem~\ref{main3} and Theorem~\ref{main4}, 
the invariant $J^{\alpha}(\sigma)$ satisfies the required 
property of Conjecture~\ref{kotae}.

\subsection*{Acknowledgement}
The author thanks Michiaki Inaba for useful discussions. He is supported 
by Japan Society for the Promotion of Sciences Research 
Fellowships for Young Scientists, No 1611452.

\subsection*{Convention}
Throughout this paper we work over $\mathbb{C}$. 
For a variety $X$, we denote by $D(X)$ the bounded 
derived category of coherent shaves on $X$. 
For a triangulated category $\tT$, its Grothendieck group 
is denoted by $K(\tT)$. When $\tT =D(X)$, we simply write it 
$K(X)$.

\section{Generalities on stability conditions}
The notion of stability conditions on triangulated categories 
was introduced in~\cite{Brs1} to give the mathematical 
framework for the Douglas's work on $\Pi$-stability~\cite{Dou1}, \cite{Dou2}. 
Here we collect some basic definitions and results 
in~\cite{Brs1}, \cite{Brs2}. 

\subsection{Stability conditions on triangulated categories}
 
\begin{defi}\label{general}\emph{
A stability condition on a triangulated category $\tT$ 
consists of data $\sigma =(Z, \pP)$, 
where $Z\colon K(\tT)\to \mathbb{C}$ is a linear map, 
and $\pP (\phi)\subset \tT$ is a full additive subcategory 
for each $\phi \in \mathbb{R}$, 
which satisfy the following:}
\emph{
 \begin{itemize}
 \item $\pP (\phi +1)=\pP (\phi)[1].$ 
\item  If $\phi _1 >\phi _2$ and $A_i \in \pP (\phi _i)$, then 
$\Hom (A_1, A_2)=0$. 
\item  If $E\in \pP (\phi)$ is non-zero,
 then $Z(E)=m(E)\exp (i\pi \phi)$ for some 
$m(E)\in \mathbb{R}_{>0}$. 
\item For a non-zero object $E\in \tT$, we have the 
following collection of triangles:
$$\xymatrix{
0=E_0 \ar[rr]  & &E_1 \ar[dl] \ar[rr] & & E_2 \ar[r]\ar[dl] & \cdots \ar[rr] & & E_n =E \ar[dl]\\
&  A_1 \ar[ul]^{[1]} & & A_2 \ar[ul]^{[1]}& & & A_n \ar[ul]^{[1]}&
}$$
such that $A_j \in \pP (\phi _j)$ with $\phi _1 > \phi _2 > \cdots >\phi _n$. 
\end{itemize}   }
\end{defi} 
We denote 
$\phi _{\sigma}^{+}(E)=\phi _1$ and $\phi _{\sigma}^{-}(E)=\phi _n$.
The non-zero objects of $\pP (\phi)$ are called \textit{semistable of 
phase} $\phi$, and
the objects $A_j$ are called \textit{semistable factors}
 of $E$ with 
respect to $\sigma$. 
For an object $E\in \tT$ the \textit{mass}
$m_{\sigma}(E)\in \mathbb{R}_{>0}$ is defined by 
$$m_{\sigma}(E)=\sum _{i=1}^n \lvert Z(A_i) \rvert.$$
The following proposition is useful in constructing stability conditions. 
\begin{prop}\emph{\bf{\cite[Proposition 4.2]{Brs1}}}\label{tstru}
Giving a stability condition on $\tT$ is equivalent to giving a heart of a
bounded t-structure $\aA \subset \tT$, and a group homomorphism 
$Z\colon K(\tT)\to \mathbb{C}$ called a stability function,
 such that for a non-zero object $E\in \aA$ one has
$$ Z(E)\in
\{ r\exp (i\pi \phi) \mid r>0, 0<\phi \le 1\},$$
and the pair $(Z, \aA)$ 
satisfies the Harder-Narasimhan 
property.
\end{prop}
For the Harder-Narasimhan property, we refer~\cite[Definition 2.3]{Brs1}. 
For a non-zero object $E\in \aA$, one can find $\phi(E) \in (0,1]$
which satisfies $Z(E) \in \mathbb{R}_{>0}e^{i\pi \phi(E)}$. We also 
call $\phi(E)$ the \textit{phase} of $E$. 
The correspondence of Proposition~\ref{tstru}
is given by 
$$(Z, \pP) \longmapsto (Z, \pP ((0,1])).$$
Here for an interval $I\subset \mathbb{R}$, the subcategory 
$\pP (I)\subset \tT$ is defined to be the smallest
extension closed subcategory 
which contains $\pP (\phi)$ for $\phi \in I$.  
In particular $\pP ((0,1])$ is a heart of a t-structure on 
$\tT$, and similarly
 $$\aA _{\phi}=\pP ((\phi -1, \phi]),$$ is also 
a heart of a t-structure for any $\phi \in \mathbb{R}$. 
(See~\cite[Section 3]{Brs1}.) 
On the other hand for $\phi _1, \phi _2\in \mathbb{R}$ with 
$0< \phi _2- \phi _1 <1$, the category 
$\pP ((\phi _1, \phi _2])$ is only a quasi-abelian 
category. We say a morphism $E_1 \to E_2$ in $\pP ((\phi _1, \phi _2])$
\textit{strict epimorphism} if it fits into 
the triangle $E_3 \to E_1 \to E_2$ with $E_3 \in \pP ((\phi _1, \phi _2])$. 
For the detail, one can 
consult~\cite[Section 4]{Brs1}, especially~\cite[Lemma 4.3]{Brs1}.

\subsection{The space of stability conditions}\label{Space}
The set of stability conditions which satisfy the technical condition
\textit{local finiteness}~\cite[Definition 5.7]{Brs1}
is denoted by $\Stab (\tT)$. 
It is shown in~\cite[Section 6]{Brs1} that $\Stab (\tT)$ has a 
natural topology. 
In fact for $\sigma \in \Stab (\tT)$
 and $\varepsilon>0$, there
 is a subset 
 \begin{align}\label{opbasis}
 B_{\varepsilon}(\sigma)\subset \Stab (\tT), \end{align}
and $\{B_{\varepsilon}(\sigma)\}_{\varepsilon, \sigma}$ gives 
an open basis of $\Stab (\tT)$. 
We refer~\cite[Section 6]{Brs1} for the construction of 
$B_{\varepsilon}(\sigma)$.  
Here we only note that for 
$\tau =(W, \qQ)\in B_{\varepsilon}(\sigma)$, one has 
$$\qQ (\phi)\subset \pP ((\phi -\varepsilon, \phi +\varepsilon)),$$
for any $\phi \in \mathbb{R}$. (See~\cite[Lemma 6.1]{Brs1}.)
Forgetting the information of $\pP$, we have the map
$$\zZ \colon \Stab (\tT) \longrightarrow \Hom _{\mathbb{Z}}(K(\tT), 
\mathbb{C}).$$
\begin{thm}\emph{\bf{\cite[Theorem 1.2]{Brs1}}}\label{lois}
For each connected component $\Sigma \subset \Stab (\tT)$, 
there exists a linear subspace $V(\Sigma )\subset \Hom _{\mathbb{Z}}(K(\tT), 
\mathbb{C})$ with a norm such that $\zZ$ restricts to a local homeomorphism, 
$\zZ \colon \Sigma \to V(\Sigma)$. 
\end{thm}
Let $\widetilde{\GL}^{+}(2, \mathbb{R})$ be the 
universal cover of $\GL ^{+}(2, \mathbb{R})$. There is 
the right action of $\widetilde{\GL}^{+}(2, \mathbb{R})$, and 
the left action of the group $\Auteq (\tT)$ on
 $\Stab (\tT)$~\cite[Lemma 8.2]{Brs1}.
 By the description in \textit{loc.cite.}, the action of
 $\widetilde{\GL}^{+}(2, \mathbb{R})$ does 
 not change the set of semistable objects. 
 The subgroup $\mathbb{C}\subset \widetilde{\GL}^{+}(2, \mathbb{R})$
 acts on $\Stab (\tT)$ faithfully. Explicitly for $\lambda \in \mathbb{C}$
 and $\sigma =(Z, \pP)$, $\lambda(\sigma)=(Z', \pP ')$ with 
 \begin{align}\label{action}
 Z'(\ast)= e^{-i\pi \lambda}Z(\ast), \quad 
 \pP '(\phi)=\pP (\phi +\Ree \lambda).\end{align}

\subsection{Numerical stability conditions}
In general $\Stab (\tT)$ is infinite dimensional. So usually 
we consider the space of \textit{numerical stability conditions}.
(See~\cite[Section 4]{Brs2}.)
Let $X$ be a smooth projective variety.
Recall that we have the pairing, 
$$\chi \colon D(X)\times D(X) \ni (E,F) \longmapsto 
\chi (E,F)=\sum _{i\in \mathbb{Z}}(-1)^i \dim \Hom (E, F[i])\in 
\mathbb{Z},$$
and it descends to the paring on $K(X)$. 
\begin{defi}\label{numGro}\emph{
We define the numerical Grothendieck group 
$\nN (X)$ to be the quotient group, 
$$\nN (X)=K(X)/\equiv, $$
where $E_1 \equiv E_2$ if and only if $\chi (E_1, F)=\chi (E_2,F)$
for any $F\in K(X)$. 
A stability condition $\sigma =(Z, \pP)$ on $D(X)$ is \textit{numerical}
if $Z\colon K(X) \to \mathbb{C}$ factors 
through
$$Z\colon K(X) \lr \nN (X) \lr \mathbb{C}.$$}
\end{defi}
The set of locally finite numerical stability conditions is 
denoted by $\Stab (X)$. There exists a map~\cite[Theorem 4.1]{Brs2},
$$\zZ \colon \Stab (X) \lr \nN (X)_{\mathbb{C}}^{\ast},$$
and since the dimension of $\nN (X)_{\mathbb{C}}$ is finite, 
 any connected 
component of $\Stab (X)$ is a complex manifold. 
In this paper, we introduce the notion of 
\textit{algebraic stability conditions}, 
whose definition is not seen in the literatures. 
\begin{defi}\emph{
We call a stability condition $\sigma =(Z, \pP)\in \Stab (X)$
algebraic if the image of $Z\colon \nN (X) \to \mathbb{C}$ is 
contained in $\mathbb{Q}\oplus  \mathbb{Q}i$. }
\end{defi}
If $\sigma =(Z, \pP)$ is algebraic, then the image of 
$Z$ is discrete and
the abelian category $\pP ((0,1])$ is 
noetherian~\cite[Proposition 5.0.1]{AP}. 
Here we put a notation and an easy remark. 
\begin{defi}\emph{ Define $\iI \subset \mathbb{R}$ to be
\begin{align}\label{rato}
\iI=\{ \phi \in \mathbb{R} 
\mid \mbox{there exists a rational point in }
 \mathbb{R}_{>0}e^{i\pi \phi} \}. \end{align}
 } \end{defi}
Note that $\iI$ is a dense countable subset in $\mathbb{R}$. 

\begin{rmk}\label{easy}\emph{
Let us take an algebraic 
stability condition $\sigma =(Z, \pP)$ and 
$\phi \in \iI$. By (\ref{action}), 
 we can find $g\in \mathbb{C}$ such that 
$g(\sigma)=(Z', \pP ')$ is also algebraic and 
$\pP '((0,1])=\pP ((\phi -1, \phi])$. 
Hence $\aA _{\phi}=\pP ((\phi -1, \phi])$ is also 
noetherian for $\phi \in \iI$. }
\end{rmk}

\subsection{Wall and chamber structures}
Let $\Stab ^{\ast}(X)$ be one of the connected components
of $\Stab (X)$. 
We use the wall and chamber structure on 
the space $\Stab ^{\ast}(X)$. For the detail one can 
consult~\cite[Section 9]{Brs2}.
For a fixed $\sigma \in \Stab ^{\ast}(X)$, 
we say a subset $\sS \subset D(X)$ has \textit{bounded mass} if there 
exists $m>0$ such that $m_{\sigma}(E)\le m$ for any $E\in \sS$. 
Note that this
 notion does not depend on a choice of $\sigma \in \Stab ^{\ast}(X)$. 
 (See~\cite[Definition 9.1]{Brs2}.)
The following is a slight generalization of~\cite[Proposition 9.3]{Brs2}. 
\begin{prop}\label{massfin}
Assume that for any bounded mass subset $\sS \subset D(X)$,
the numerical classes 
\begin{align}\label{katakuna}\{ [E] \in \nN (X) \mid E\in \sS \},\end{align}
is a finite set. Then for any compact subset $\mathfrak{B}\subset 
\Stab ^{\ast}(X)$, there exists a finite number of real codimension one 
submanifolds $\{\wW _{\gamma} \mid \gamma \in \Gamma \}$
 on $\Stab ^{\ast}(X)$ such that 
 if $\Gamma '$ is a subset of $\Gamma$ and $\cC$ is one of the 
 connected components, 
 \begin{align}\label{component}
 \cC \subset \bigcap _{\gamma \in \Gamma '}(\mathfrak{B}\cap 
 \wW _{\gamma}) \setminus 
 \bigcup _{\gamma \notin \Gamma '}\wW _{\gamma},\end{align}
 then if $E\in \sS$
  is semistable in some $\sigma \in \cC$, then it is semistable 
 for all $\sigma \in \cC$. 
\end{prop}
\begin{proof}
The statement is not seen in the literatures. However
the proof is a straightforward adaptation of the 
proof of~\cite[Proposition 9.3]{Brs2} and we leave the readers to check
the detail. 
Here we only recall the construction of the walls $\{ \wW _{\gamma}\}_{\gamma \in \Gamma}$, since it will be needed later. 
 For a bounded mass subset $\sS \subset D(X)$, 
Bridgeland~\cite[Proposition 9.3]{Brs2} considered another 
bounded mass subset $\sS \subset \sS ' \subset D(X)$, 
$$\sS '= \{ A\in D(X) \mid \mbox{there is some }\sigma \in \mathfrak{B}
\mbox{ and }E\in \sS \mbox{ such that }m_{\sigma}(A)\le m_{\sigma}(E)
\},$$
and let 
$v_1, \cdots, v_n \in \nN (X)$ be the numerical classes of $\sS'$. 
Let $\Gamma$ be the set of pairs $(v_i, v_j)$ such that $v_i$ and 
$v_j$ are not proportional in $\nN (X)$. Then for $\gamma =(v_i, v_j)\in 
\Gamma$, $\wW _{\gamma}$ is defined to be 
\begin{align}\label{wall}
\wW _{\gamma}=\{ \sigma =(Z, \pP)\in \Stab ^{\ast}(X) \mid 
Z(v_1)/Z(v_2) \in \mathbb{R}_{>0}\}.\end{align}
\end{proof}
It is proved in~\cite[Lemma 9.2]{Brs2} that
the assumption of Proposition~\ref{massfin} is satisfied 
when $X$ is a K3 surface or an abelian surface.  

We say a connected component $\Stab ^{\ast}(X)$ is
\textit{full} if the image of the map 
$\Stab ^{\ast}(X) \to \nN (X)_{\mathbb{C}}^{\ast}$ is an open  
subset of $\nN (X)_{\mathbb{C}}^{\ast}$.
Note that if $\Stab ^{\ast}(X)$ is full, then 
the subset of algebraic stability conditions is 
dense in $\Stab ^{\ast}(X)$.
Here we give the following easy lemma. 
\begin{lem}\label{put}
Assume that $\Stab ^{\ast}(X)$ is full. 
Let $\mathfrak{B}^{\circ}$ be an open subset of $\Stab ^{\ast}(X)$ 
and its closure $\mathfrak{B}$ is compact. Then for a
connected component $\cC$ of (\ref{component}), 
the set of points $\sigma \in \cC$ which are algebraic is 
dense in $\cC$. 
\end{lem}
\begin{proof}
There is no proof in the literatures, however
by the description of the walls (\ref{wall}), 
 it is easy to check that any intersection 
$\cap _{\gamma \in \Gamma '}\wW _{\gamma}$ contains a dense 
subset of algebraic stability conditions.  
\end{proof}

\section{Moduli stacks of semistable objects}\label{stack}
The purpose of this section is to establish the general 
arguments to study the moduli stacks of semistable objects. 
Throughout this section, $X$ is a smooth projective variety over 
$\mathbb{C}$, 
and $S$ is a $\mathbb{C}$-scheme. 
We always assume $S$ is connected. 
For an object $\eE \in D(X\times S)$
and a $S$-scheme $T\to S$, we denote by $\eE _T$ the derived pull-back of 
$\eE$
to $X\times T$. We denote 
$$p\colon X\times S \to X, \quad q_S \colon 
X\times S \to S,$$
the projections respectively. 
For a set of objects $\sS \subset D(X)$,
we say it is \textit{bounded} if there is a $\mathbb{C}$-scheme $Q$ 
of finite type and $\fF \in D(X\times Q)$ such that any object $E\in \sS$
is isomorphic to $\fF _{q}$ for some $q\in Q$. Also we say a map 
$$\nu \colon \sS \lr \mathbb{R},$$
is bounded (resp bounded above, bounded below) 
if there is $c\in \mathbb{R}$ such that 
$\lvert \nu (E) \rvert \le c$. (resp $\nu (E)\le c$, $\nu (E)\ge c$.)
For the generalities of Artin stacks, one can consult~\cite{GL}. 
In this section, we work over 
a connected component $\Stab ^{\ast}(X)\subset \Stab (X)$,
which satisfies the following assumption. 
\begin{assum}\label{nbd}
\emph{
\begin{itemize}
\item For any bounded mass subset $\sS \subset D(X)$, the set 
of numerical classes (\ref{katakuna}) is finite. 
\item 
There is a subset $\vV \subset \Stab ^{\ast}(X)$
which consists of algebraic stability conditions
and satisfies the following:
for any algebraic
$\sigma \in \Stab ^{\ast}(X)$, 
there exist $\Phi \in \Auteq D(X)$ and $g\in \widetilde{\GL}^{+}(2,\mathbb{R})$
such that $g\circ \Phi (\sigma)$ is also algebraic and
contained in $\overline{\vV}$. \end{itemize}}
\end{assum}
The above assumption is known to hold in several examples.
For instance, if $X$ is an elliptic curve, one can take 
$\vV$ to be just one point of an algebraic stability condition~\cite{Brs1}.  
When $X$ is a K3 surface or an abelian surface, we
will see in the next section that Assumption~\ref{nbd} is 
satisfied.

\subsection{Openness of stability conditions}
Let $\mM$ be the 2-functor 
$$\mM \colon (\Sch /\mathbb{C}) \lr (\groupoid),$$
which sends a $\mathbb{C}$-scheme $S$ to the groupoid
$\mM (S)$ whose objects consist of $\eE \in D(X\times S)$ which 
is relatively perfect~\cite[Definition 2.1.1]{LIE} and
satisfies 
\begin{align}\label{perf}
\Ext ^{i}(\eE _s,\eE _s)=0, \mbox{ for all } i<0
\mbox{ and }s\in S.\end{align}
For the detail we refer~\cite{LIE}. 
Lieblich showed the following. 
\begin{thm}\emph{\bf{\cite{LIE}}}\label{artin}
The 2-functor $\mM$ is an Artin stack of locally finite 
type over $\mathbb{C}$. 
\end{thm}
Let us fix $\sigma =(Z, \pP) \in \Stab ^{\ast} (X)$, $\phi \in \mathbb{R}$ and 
$\alpha \in \nN (X)$. 
Note that any object $E\in \pP (\phi)$ satisfies (\ref{perf}). Thus 
it is possible to
define the following. 
\begin{defi} \emph{
We define $M^{(\alpha, \phi)}(\sigma)$ to be the set of 
$\sigma$-semistable objects of phase $\phi$ and numerical type $\alpha$, and
$$\mM _{}^{(\alpha, \phi)}(\sigma) \subset \mM,$$
the substack of objects in $M^{(\alpha, \phi)}(\sigma)$. }
\end{defi}
We have the following. 
\begin{lem}\label{artinstack}
Assume $M^{(\alpha, \phi)}(\sigma)$ is bounded and 
$\mM ^{(\alpha, \phi)}(\sigma)$ is an open substack of 
$\mM$. Then $\mM ^{(\alpha, \phi)}(\sigma)$ is an Artin stack 
of finite type over $\mathbb{C}$. 
\end{lem}
\begin{proof}
Let $M \to \mM$ be an atlas of $\mM$. The openness of 
$\mM ^{(\alpha, \phi)}(\sigma)$ implies there is an open subset 
$M^{\circ}\subset M$ which gives a surjective smooth morphism 
$M^{\circ} \to \mM ^{(\alpha, \phi)}(\sigma)$. Furthermore
the boundedness of $M^{(\alpha, \phi)}(\sigma)$ implies there 
is a surjection $M' \to M^{\circ}$ from a finite type $\mathbb{C}$-scheme 
$M'$. This implies $M^{\circ}$ is also of finite type, and it gives
an atlas of $\mM ^{(\alpha, \phi)}(\sigma)$. 
\end{proof}

Our purpose here is to 
give the sufficient condition for $\mM ^{(\alpha, \phi)}(\sigma)$ to 
be an open substack of $\mM$. 
We consider the following claim. 
\begin{claim}\label{op2}
For a smooth
 quasi-projective variety $S$ and
 $\eE \in \mM (S)$, 
 assume that the 
locus 
\begin{align}\label{locus}
S^{\circ}=
\{ s\in S \mid \eE _s \mbox{ is of numerical type }\alpha \mbox{ and }
\eE _s \in \pP (\phi) \},\end{align}
is not empty. Then there is an open subset $U\subset S$ which is 
contained in 
$S^{\circ}$.
\end{claim}

By the following lemma, it is enough to consider Claim~\ref{op2}. 

\begin{lem}\label{opp}
Assume Claim~\ref{op2} is true.
 Then $\mM ^{(\alpha, \phi)}(\sigma)$ is 
an open substack of $\mM$. 
\end{lem}
\begin{proof}  
By Theorem~\ref{artin}, 
it suffices to show that
 for an arbitrary affine $\mathbb{C}$-scheme $S$ of finite type 
and $\eE \in \mM(S)$, the locus 
(\ref{locus})  
is open in $S$. 
Assume Claim~\ref{op2} is true and take an affine $\mathbb{C}$-scheme 
$S$ of finite type and $\eE \in \mM (S)$. Assume that the locus (\ref{locus}) 
 is not empty. 
Let $g\colon S' \to S$ be a resolution of singularities. 
Note that the locus $S^{'\circ}\subset S'$ determined by
 $\eE _{S'}\in \mM (S')$ and (\ref{locus}) 
 is not empty because $g$ is surjective. 
Applying 
Claim~\ref{op2} to $\eE_{S'}$,
 there is an open subset $U _1 '\subset S'$ such 
that $U_1 '\subset S^{'\circ}$. Restricting to the locus where 
$g$ is an isomorphism, we obtain an open subset $U_1 \subset S$ such that 
$U_1 \subset S^{\circ}$.  
Let $Z_1 =S\setminus U_1$. If $Z_1 \cap 
S^{\circ}$ is empty, we have $S^{\circ}=U_1$. Otherwise 
take the pull-back $\eE_{Z_1} \in \mM (Z_1)$ and apply
the same argument. Then we
 obtain an open subset $U_2 \subset Z_1 \cap S^{\circ}$ 
in $Z_1$ and a closed subset $Z_2=Z_1 \setminus U_2$, which is also 
closed in $S$. 
 Repeating this argument, we get a sequence of closed subsets in $S$, 
 $$ \cdots  \subset Z_n \subset Z_{n-1} \subset \cdots \subset Z_1, $$
 which must be terminate because $S$ is noetherian. 
 Then $Z=\cap _i Z_i$ is a closed subset of $S$ and we have 
 $S^{\circ}=S\setminus Z$. Therefore $S^{\circ}$ is open. 
 \end{proof}

\subsection{Sheaf of t-structures}
Here we introduce the sheaf of t-structures studied by 
D.Abramovich and A.Polishchunk~\cite{AP}.
Let $\aA \subset D(X)$ be a heart of a bounded t-structure 
and assume that $\aA$ is noetherian. 
 Take 
a smooth projective 
variety $S$ and an ample line bundle $\lL \in \Pic (S)$. 

\begin{thm}\emph{\bf{\cite[Theorem 2.6.1]{AP}}} 
The subcategory 
$$
\aA _{S} =\{ F\in D(X\times S) \mid 
\dR p_{\ast}(F\otimes \lL ^{n}) \in \aA \mbox{ for all }
n\gg 0 \},$$
is a
 heart of a bounded t-structure on $D(X\times S)$,
 independent of a choice of $\lL$. Furthermore it is
 a noetherian abelian category. 
 \end{thm}
 The subcategory $\aA _S \subset D(X\times S)$ extends to 
 a sheaf of bounded t-structures~\cite[Theorem 2.7.2]{AP}, 
 i.e. for an open subset $j\colon U \subset S$, there exists a heart of 
 a bounded t-structure $\aA _U \subset D(X\times U)$ such that 
 $$(\id \times j)^{\ast} \colon D(X\times S) \lr D(X\times U),$$
 takes $\aA _S$ to $\aA _U$. Moreover it is shown in~\cite[Lemma 3.2.1]{AP}
 that $\aA _S \to \aA _U$ is essentially surjective, and $\aA _U$ does 
 not depend on a projective compactification $U\subset S$.  
 Thus one can define $\aA _S$ for a smooth quasi-projective 
 variety $S$. 
 One of the necessary fact for our purpose is the following 
 \textit{open heart property}. 
 \begin{thm}\emph{\bf{\cite[Theorem 3.3.2]{AP}}}\label{heart}
 For a smooth quasi-projective variety $S$ and 
$\eE\in D(X\times S)$, assume there exists $s\in S$ such that 
 $\eE _s \in \aA$. Then there exists an open neighborhood $s\in U\subset S$
 such that $\eE _U \in \aA _U$. 
 \end{thm}
 We say $\eE \in \aA _S$ is \textit{t-flat} if for any $s\in S$
 one has $\eE _s \in \aA$. Since $U\mapsto \aA _U$ is a sheaf of 
 t-structures, if $\eE \in \mM (S)$ satisfies $\eE _s \in \aA$
 for all $s\in S$, then Theorem~\ref{heart} and~\cite[Lemma 2.1.1]{AP} 
 show that $\eE \in \aA _S$ and it is t-flat. 
 For a closed point $s\in S$ and the inclusion $i_s \colon X\times \{ s \}
 \hookrightarrow X\times S$, 
 it is shown
 in~\cite[Lemma 2.5.3]{AP} that 
 $$\dL i_{s} ^{\ast} \colon D(X\times S) \to D(X),$$
 is right t-exact with respect to the t-structures with hearts 
 $\aA _S$, $\aA$ respectively. 
 Thus one has the following lemma. 
 \begin{lem}\label{shoex}
 Let $0\to \hH \to \eE \to \fF \to 0$ be an exact sequence in 
 $\aA _S$ and assume that 
 $\eE$, $\fF$
 are t-flat. Then $\hH$ is also t-flat. 
 \end{lem}
 For our purpose, we have to consider the following problem
 called \textit{generic flatness problem}.
 \begin{prob}\emph{\bf{\cite[Problem 3.5.1]{AP}}}\label{pro}
 For $\eE\in \aA _S$, is there an open subset $U\subset S$ such 
 that for each $s\in U$, we have $\eE_s \in \aA$? 
 \end{prob}
 
 \begin{rmk}\label{same}
 \emph{If Problem~\ref{pro} is true, the same argument of Lemma~\ref{opp} 
 shows the following: for an arbitrary} 
 $\mathbb{C}$\emph{-scheme} $S$ \emph{of finite type, the points}
 $s\in S$ \emph{on which }$\eE _s \in \aA$
 \emph{is in fact open.}
 \end{rmk}
 
 In~\cite{AP}, there is a partial result for Problem~\ref{pro}. 
 
 \begin{prop}\emph{\bf{\cite[Proposition 3.5.3]{AP}}}\label{dense}
 For $\eE\in \aA _S$, there is a dense subset $U\subset S$ such 
 that for each $s\in U$, we have $\eE_s \in \aA$.
 \end{prop}

 The generic flatness problem requires $U$ to be 
 open in Zariski topology. 
 Let us take an algebraic stability condition 
 $\sigma =(Z, \pP) \in \Stab ^{\ast} (X)$.   
The purpose here is to reduce Claim~\ref{op2} to Problem~\ref{pro}.

\begin{lem}\label{op3} 
Let $S$ be a smooth quasi-projective variety, $\phi \in \mathbb{R}$
and $\alpha \in \nN (X)$.  

(i) For $\eE \in \mM (S)$, 
assume the locus $S^{\circ}$ defined by (\ref{locus}) is non-empty. 
Then $S^{\circ}$ is 
dense in $S$. 

(ii) In the same situation of (i), assume Problem~\ref{pro} is true 
 for $\aA _{\phi}=\pP((\phi -1, \phi])$. 
 Then $S^{\circ}$ contains an open subset of $S$. 
(Thus Claim~\ref{op2} is true for this $\phi \in \mathbb{R}$.)  
\end{lem}
\begin{proof}  
(i) 
 Because $S^{\circ}$ is non-empty and $\sigma$ is algebraic,
  we have $\phi \in \iI$. Hence $\aA _{\phi}= \pP ((\phi -1, \phi])$ 
  is noetherian by Remark~\ref{easy}. 
 Let us take $s\in S^{\circ}$.  
 Note that
$\eE _s \in \pP (\phi)\subset \aA _{\phi}$. Thus by Lemma~\ref{easy} and
Theorem~\ref{heart},
 there exists 
an open subset $s\in U\subset S$ such that $\eE _U \in \aA _{\phi,U}$. 
Therefore by Proposition~\ref{dense}
there exists dense subset $U'\subset U$ such that for 
$s'\in U'$, we have $\eE _{s'}\in \aA _{\phi}$.
 Since $\eE _{s'}$ is numerically 
equivalent to $\eE _s$, 
we have $Z(\eE _{s'})\in \mathbb{R}_{>0}e^{i\pi \phi}$. 
This implies $\eE _{s'} \in \pP (\phi)$, hence $U' \subset S^{\circ}$.

(ii) If we assume the generic flatness for
$\aA _{\phi}$, then we can take $U'$ in the 
proof of (i) to be open.  
\end{proof}

\subsection{Boundedness of semistable objects}
Here we discuss the boundedness of semistable objects and 
certain quotient objects. 
We fix an algebraic stability condition
$\sigma =(Z, \pP) \in \Stab ^{\ast} (X)$, and consider the following problem. 

\begin{problem}\label{further}
Is the 
set of objects
$M^{(\alpha,\phi)}(\sigma)$
bounded, for any $\alpha \in \nN (X)$ and $\phi \in \mathbb{R}$ ?
\end{problem}
Let $\aA =\pP ((0,1])$. We show the following. 
\begin{lem}\label{quotb}
Assume Problem~\ref{further} is true for a fixed $\sigma$. 
Then for any $\phi \in (0,1)$ and $G\in \aA$,
 the following set of objects, 
$$Q(G,\phi) =\{ 
E\in \aA  \mid \mbox{there exists a surjection }
G \twoheadrightarrow E \mbox{ in } \aA  \mbox{ and }
\phi (E) \le \phi \},$$
is bounded.

\end{lem}
\begin{proof} 
For $E\in 
Q(G, \phi)$, let $F_1, F_2, \cdots, F_{n(E)}$ be the semistable 
factors of $E$ in $\sigma$ such that $F_i \in \pP (\phi _i)$ and 
$\phi _1 > \phi _2 > \cdots > \phi _{n(E)}$. We have 
\begin{align}\label{night}
\sum _{i=1}^{n(E)} \Imm Z(F_i) = \Imm Z(E) \le \Imm Z(G).
\end{align}
Note that $\Imm Z(F_i) >0$ except $i=1$. Because $\sigma$ is algebraic, 
the image $\nN (X) \stackrel{Z}{\to} \mathbb{C} \stackrel{\Imm}{\to} 
\mathbb{R}$ is discrete. Thus (\ref{night}) implies that the map 
$E \mapsto n(E)$ on $Q(G, \phi)$ is bounded, and the following set 
\begin{align}\label{night2}
\{ \Imm Z(F_i) \in \mathbb{Q} \mid 1\le i\le n(E), E\in Q(G,\phi) \},
\end{align}
is a finite set. 

Next there exist surjections, $G \twoheadrightarrow E \twoheadrightarrow 
F_{n(E)}$ in $\aA$, so we have $\phi _{\sigma}^{-}(G) \le 
\phi _{n(E)} \le \phi _i$ for $1\le i\le n(E)$. (See~\cite[Lemma 3.4]{Brs1}.)
Thus the map on $Q(G, \phi)$,
\begin{align}\label{night3}
E \longmapsto \max \{ \Ree Z(F _i) \mid 1\le i\le n(E) \} \in \mathbb{Q},
\end{align}
 is bounded above. On the other hand since $\phi (E) \le \phi  <1$
 and
 $\Imm Z(E) \le \Imm Z(G)$, the following map on $Q(G,\phi)$,
 \begin{align}\label{night4}
 E \longmapsto \Ree Z(E) =\sum _{i=1}^{n(E)} \Ree Z(F_i),
 \end{align}
 is bounded below. 
 Combined with the fact that (\ref{night3}) is bounded above,
  the following set 
 \begin{align}\label{night5}
 \{ \Ree Z(F_i) \in \mathbb{Q} 
 \mid 1 \le i\le n(E), E\in Q(G, \phi) \},
 \end{align}
 is a finite set. 
 Then the finiteness of (\ref{night2}),
 (\ref{night5}) and Assumption~\ref{nbd} imply that the following set, 
 \begin{align}\label{night6}
 \{ [F_i] \in \nN (X) \mid 1 \le i\le n(E), E\in Q(G, \phi) \},
 \end{align}
 is a finite set. 
 Since we assume that Problem~\ref{further} is true, 
  the finiteness of (\ref{night6})
  implies that the set of objects 
 $$\{ F_i \mid 1\le i\le n(E), E\in Q(G, \phi) \},$$
 is bounded. Thus $Q(G, \phi)$ is also bounded by Lemma~\ref{useasy} below.  
\end{proof}
Here we have used the following easy lemma. 
\begin{lem}\label{useasy}
Let $\sS _i \subset D(X)$ be the sets of objects for $1\le i\le 3$ and 
 $\sS _1$, $\sS _2$ are bounded. Assume that for any object $E _3 \in \sS _3$, 
 there is $E_i \in \sS _i$ for $i=1,2$ and a triangle, 
 $$E_1 \lr E_3 \lr E_2.$$
 Then $\sS _3$ is also bounded. 
\end{lem}
The proof is easy and leave it to the reader. 
In fact it is enough to notice that 
$\Ext ^1 (E_2, E_1)$ is finite dimensional.  

Assuming Problem~\ref{pro} and Problem~\ref{further}, we can 
construct certain schemes which parameterize quotient objects. 
Let 
$\eE \in \aA _S$ be a t-flat family and take $\phi \in (0,1)$. 
 We consider the following functors, 
$$ \mathrsfs{Q} uot (\eE, \phi), \mathrsfs{S} ub (\eE, \phi)
 \colon (\Sch /S) \lr 
(\Set),
$$
such that $\mathrsfs{Q} uot (\eE , \phi)$ (resp $\mathrsfs{S} ub (\eE, \phi))$
 takes a $S$-scheme $T$ to the isomorphism classes of objects 
$\fF \in \mM (T)$ together 
with a morphism $\eE _T \to \fF$, (resp $\fF \to \eE _T$)
such that 
\begin{itemize}
\item For each closed point $t\in T$, $\fF _t$ is contained in 
$\aA$ and $\phi (\fF _t) \le \phi$. (resp $\phi (\fF _t) >
\phi$.) 
\item For each closed point $t\in T$, the induced morphism 
$\eE _t \to \fF _t$ is a surjection in $\aA$. 
(resp $\fF _t \to \eE _t$ is an injection in $\aA$.)
\end{itemize}
We show the following. 

\begin{prop}\label{quot}
For a fixed $\sigma$, assume
 Problem~\ref{pro} for $\aA =\pP ((0,1])$ and Problem~\ref{further} are
 true. 
Then for any $\phi \in (0,1)$ there
exist $S$-schemes $\qQ(\eE,\phi)$, $\sS(\eE, \phi)$
 which are of finite type 
over $S$, and morphisms over $S$, 
\begin{align*}
\qQ(\eE,\phi) &\lr \mathrsfs{Q} uot (\eE, \phi),\\
\sS(\eE, \phi) &\lr \mathrsfs{S} ub (\eE, \phi),
\end{align*}
which are surjective on $\mathbb{C}$-valued points of $\mathrsfs{Q} uot (\eE,\phi)$ 
and $\mathrsfs{S} ub (\eE, \phi)$. 
\end{prop}

\begin{proof}
First let us construct $\qQ(\eE, \phi)$. 
By~\cite[Lemma 2.6.2]{AP}, there exists an object 
$G\in \aA$, $n\in \mathbb{Z}$ and a surjection 
$G_{S}\otimes \lL ^{-n} \twoheadrightarrow \eE$ in $\aA _S$. Note that 
the induced morphism $G \to \eE _s$ is a surjection 
by Lemma~\ref{shoex}. 
By the assumption and Lemma~\ref{quotb}, 
 there is a $\mathbb{C}$-scheme $Q_1$ of finite 
type over $\mathbb{C}$ and an object $\fF \in \mM (Q_1)$ such that 
any object in $Q(G,\phi)$ is isomorphic to $\fF _{q}$ for some $q\in Q_1$. 
Let $Q_1 ^{\circ}$ be 
$$Q_1 ^{\circ}= \{ q\in Q_1 \mid \fF _{q}\in \aA \}.$$
Since we assume the generic flatness for $\aA$, the locus 
$Q_1 ^{\circ}$ is open in $Q_1$. Set $Q_2 =Q_1 ^{\circ}\times S$
and we regard it as a $S$-scheme via the projection $Q_2 \to S$. 
By~\cite[Proposition 2.2.3]{LIE}, there exists an affine
open subset $U\subset Q_2$
such that the functor $\Coh (U) \to \Coh (U)$ sending $M$ to 
\begin{align}\label{mono}
M \longmapsto \hH ^0 (\dR q_{U \ast}\dR \hH om (\eE _U, \fF _{U} \otimes 
q_{U}^{\ast}M)),\end{align}
has the form $\hH om (\widetilde{\eE}_U,M)$ for some locally free 
sheaf $\widetilde{\eE}_U$ on $U$. 
Here $\fF _U$ is the pull-back of $\fF$ via
$$U\subset Q_2 \to Q_1 ^{\circ} \subset Q_1.$$
Set 
$Q_2 ' =(Q_2 \setminus U) \coprod U$ and apply the same procedure 
to $\eE _{Q_2 '}$ and $\fF _{Q_2 '}$ repeatedly. Then we 
obtain an affine scheme of finite type $Q_3$ with a morphism
$Q_3  \to Q_2$, which is bijective on closed points, and a 
locally free sheaf $\widetilde{\eE}$ on $Q_3$ such that the functor 
$\Coh (Q_3) \to \Coh (Q_3)$ given in the same way as (\ref{mono}) has the 
form $\hH om (\widetilde{\eE}, \ast)$. Furthermore 
the functor 
$$(T\to Q_3) \longmapsto \hH ^0 (\dR q_{T\ast}\dR \hH om (\eE _T, \fF _T))
\in \Coh (T),$$
is represented by $\mathbb{V}(\widetilde{\eE})$
by~\cite[Proposition 2.2.3]{LIE}.
 Thus there exists a universal 
morphism
$\eE _{\mathbb{V}(\widetilde{\eE})} \to \fF _{\mathbb{V}(\widetilde{\eE})}$.
 Let $\hH$ be its cone, i.e. $\hH$ fits 
into the distinguished triangle in $D(X\times \mathbb{V}(\widetilde{\eE}))$, 
$$\hH \lr \eE _{\mathbb{V}(\widetilde{\eE})} \lr \fF _{\mathbb{V}(\widetilde{\eE})}.$$
For $q\in \mathbb{V}(\widetilde{\eE})$, 
note that $\fF _{q}$ is contained in $\aA$. 
Thus the induced morphism $\eE _q \to \fF _q$ is surjective in 
$\aA$ if and only if $\hH _q \in \aA$. 
Then define $\qQ(\eE,\phi)$ to be the locus, 
$$\qQ(\eE,\phi)
\cneq \{ q\in \mathbb{V}(\widetilde{\eE}) \mid \hH _q \in \aA \}.$$
Again $\qQ(\eE,\phi)$
 is an open subscheme of $\mathbb{V}(\widetilde{\eE})$, in 
 particular it is of finite 
type over $S$. The 
restriction of $\eE _{\mathbb{V}(\widetilde{\eE})} \to 
\fF _{\mathbb{V}(\widetilde{\eE})}$
to $\qQ(\eE, \phi)$ 
induces a morphism, 
$$\qQ(\eE,\phi)
 \lr \mathrsfs{Q} uot (\eE,\phi),$$ which is surjective 
on $\mathbb{C}$-valued points by the construction. 

Next we construct $\sS(\eE, \phi)$. Since $\phi (\eE _s)$
does not depend on $s\in S$, we can easily see the following:
there exists $\phi ' \in (0,1)$ such that for any $s\in S$ and 
a subobject $\hH \subset \eE _s$ in $\aA$ with $\phi (\hH)> \phi$, 
we have $\phi (\eE _s/ \hH) \le \phi '$. Let us consider 
$\qQ(\eE, \phi ')$ and the universal quotient 
$\eE _{\qQ(\eE, \phi ')} \to \fF$ on $X\times 
\qQ(\eE, \phi ')$. We consider the distinguished triangle, 
$$\hH \lr \eE _{\qQ(\eE, \phi ')} \lr \fF.$$
Note that $\hH _q \in\aA$, thus one can 
define its phase $\phi (\hH _q)\in (0,1]$. 
 Then we construct $\sS(\eE, \phi)$ 
as follows, 
$$\sS(\eE, \phi) =\{ q\in \qQ(\eE, \phi ') \mid 
\phi (\hH _{q}) > \phi \}.$$
Since $q \mapsto \phi (\hH _{q})$ is locally constant on $\qQ(\eE, \phi ')$, 
the locus $\sS(\eE, \phi)$ is a union of the 
connected components of $\qQ(\eE, \phi ')$, in particular 
of finite type over $S$. 
The induced morphism $\hH _{\sS(\eE, \phi)} \to 
\eE _{\sS(\eE, \phi)}$ gives a morphism 
$\sS(\eE, \phi) \to \mathrsfs{S} ub (\eE, \phi)$, which is  
surjective on $\mathbb{C}$-valued points. 
\end{proof}

\subsection{Generic flatness for $\aA _{\phi}=\pP ((\phi -1, \phi])$}
Again we fix an algebraic stability condition
$\sigma =(Z, \pP)\in \Stab ^{\ast} (X)$.
Here we study the generic flatness for $\pP ((\phi -1, \phi])$
under several assumptions. The purpose here is the following. 
\begin{prop}\label{purpose}
Under the same assumption as in Proposition~\ref{quot}, 
let us take $\phi \in \iI$. 
Then Problem~\ref{pro} is true for
$\aA _{\phi}=\pP ((\phi -1, \phi])$. 
\end{prop}
\begin{proof}
For $\eE \in \aA _{\phi,S}$, let us find an open subset $U\subset S$
on which $\eE _s\in \aA _{\phi}$. 
We may assume $0<\phi \le 1$. 
By~\cite[Lemma 3.2.1]{AP} we may also assume $S$ is projective,
and let $\lL \in \Pic (S)$ be an ample line bundle.  
Then $\eE \in \aA _{\phi,S}$
implies, 
\begin{align}\label{kimi}
\dR p_{\ast}(\eE \otimes \lL ^n)\in \pP ((\phi -1, \phi]),\end{align}
for $n\gg 0$. 
For $\aA =\pP ((0,1])$,  
we denote by 
$H_{\aA}^i(\ast)$, $H_{\aA _S}^i (\ast)$ the $i$-th cohomology functors
on $D(X)$, $D(X\times S)$
with respect to the t-structures with hearts $\aA$, $\aA _S$ respectively.
Then (\ref{kimi}) implies 
\begin{align}\label{love}
H_{\aA}^i (\dR p_{\ast}(\eE \otimes \lL ^n))=0 \mbox{ unless } 
i=0,1.\end{align}
On the other hand, we have 
\begin{align}\label{stop}
\dR p_{\ast}(H^i _{\aA _S}(\eE \otimes \lL ^n))=
\dR p_{\ast}(H^i _{\aA _S}(\eE)\otimes \lL ^n) \in \aA,\end{align}
for $n\gg 0$. The first equality comes from~\cite[Proposition 2.1.3]{
AP}.
Thus (\ref{love}) and (\ref{stop}) imply 
\begin{align}\label{kataki}\dR p_{\ast}(H^i _{\aA _S}(\eE)\otimes \lL ^n) =0
\mbox{ unless } i=0,1,\end{align}
for $n\gg 0$. 
It is easy to deduce from (\ref{kataki}) that $H^i _{\aA _S}(\eE) =0$
unless $i=0,1$, by using the standard t-structure on $D(X\times S)$. 
For $i=0,1$, denote $\eE ^i =H^i _{\aA _S}(\eE)\in \aA _S$. Since we assume 
the generic flatness for $\aA$, there exists an open set $S_1 \subset S$
such that for $s\in S_1$ one has $\eE ^i _s \in \aA$. 
Since we have the distinguished triangle $\eE ^0 \to \eE \to \eE ^1[-1]$, 
we have the distinguished triangle in $D(X)$, 
$$\eE ^0 _s \lr \eE _s \lr \eE ^1 _s [-1],$$
for $s\in S_1$.
Hence for $s\in S_1$, $\eE _s \in \aA _{\phi}$ is equivalent to the following, 
\begin{align}\label{phant}
\eE _s ^0 \in \pP ((0,\phi]), \quad \eE _s ^1 \in \pP ((\phi, 1]).
\end{align}
Thus it is enough to find an open set $U\subset S_1$ where 
(\ref{phant}) holds. 
Note that by Proposition~\ref{dense}, the set of points $s\in S_1$
on which $(\ref{phant})$ hold is dense in $S$. 
First let us consider the locus where 
$\eE _s ^1 \in \pP ((\phi, 1])$ holds.
Let $\pi _{\phi}$ be the composition of the morphisms, 
$$
\pi _{\phi } \colon \qQ(\eE ^1, \phi) \lr 
\mathrsfs{Q} uot (\eE ^1, \phi) \lr
S,$$
constructed in Proposition~\ref{quot}.
Note that $\eE _s ^1 \in \pP ((\phi, 1])$ if and only if 
there is no surjection $\eE _s ^1 \twoheadrightarrow F$ in 
$\aA$ with $\phi (F) \le \phi$. 
Thus $\eE _s ^1 
\in \pP ((\phi, 1])$ if and only if $s \notin \im \pi _{\phi}$, 
and such points are dense in $S_1$. 
This implies  
$\pi _{\phi}$ is not dominant. Because $\qQ(\eE ^1, \phi)$ is 
of finite type, there is an open subset 
$U\subset S_1 \setminus \im \pi _{\phi}$, and (\ref{phant}) holds on $U$. 

We can argue in a similar way (using $\sS (\eE, \phi)$ instead 
of $\qQ(\eE, \phi)$)
to find an open subset $U\subset S_1$
where $\eE _s ^0 \in \pP ((0,\phi])$ holds.
 We leave the detail to the reader. 
\end{proof}

\subsection{Sufficient conditions for $\mM ^{(\alpha, \phi)}(\sigma)$
 to be an Artin stack of finite type}
 Here we give the sufficient condition for $\mM ^{(\alpha, \phi)}(\sigma)$ to 
 be an Artin stack of finite type. First let us consider a slight 
 generalization of Proposition~\ref{quot}. 
As before we fix an algebraic stability condition
 $\sigma =(Z, \pP)\in \Stab ^{\ast} (X)$, 
and $\aA =\pP ((0,1])$.  
Take $\phi _0, \phi _1 \in \iI$ with $\phi _1 -\phi _0 <1$, 
and
$\eE \in \mM (S)$ which satisfies $\eE _s \in \pP ((\phi _0, \phi _1])$
 for all $s\in S$. We
define the functor,
 $$\mathrsfs{Q} uot (\eE, \phi _0, \phi _1)\colon 
(\Sch /S) \lr (\Set),$$
 by associating a $S$-scheme $T$ to the set of 
isomorphism classes of $\fF \in \mM(T)$ together with a morphism
$\eE _T \to \fF $ such that for each $t\in T$, the induced morphism 
$\eE _t \to \fF _t$ is a strict epimorphism in $\pP ((\phi _0, \phi _1])$. 
We need the following. 

\begin{lem}\label{ad}
Under the same assumption as in Proposition~\ref{quot}, 
take $\phi _0, \phi _1 \in \iI$ as above. 
Then there exists a $S$-scheme $\qQ(\eE, \phi _0, \phi _1)$
of finite type over $S$
and a morphism 
$$\qQ(\eE, \phi _0, \phi _1) \lr \mathrsfs{Q} uot (\eE, \phi _0, \phi _1),$$
which is surjective on $\mathbb{C}$-valued points. 
 \end{lem}
 \begin{proof}
 Let us take $\phi _2 \in \iI$ which satisfies 
 $\phi _0, \phi _1 \in (\phi _2 -1, \phi _2)$, and 
 $G\in \aA _{\phi _2}= \pP ((\phi _2 -1, \phi _2])$. 
 By Remark~\ref{easy}, one can apply Lemma~\ref{quotb} and conclude that 
 the following set of objects 
 $$Q(G, \phi _0, \phi _1)=\{ E\in \pP ((\phi _0, \phi _1]) \mid 
\mbox{there exists a surjection } G \twoheadrightarrow E 
\mbox{ in } \aA _{\phi _2} \}, $$
is bounded. 
 As in Proposition~\ref{quot}, there exists a surjection 
$G_S \otimes \lL ^{-n} \twoheadrightarrow \eE$ in $\aA _{\phi _2,S}$ for 
some $n\in \mathbb{Z}$ and $\lL \in \Pic (S)$ is an ample line bundle.  
By the boundedness of $Q(G, \phi _0, \phi _1)$, there 
exists a $\mathbb{C}$-scheme $Q_1$ of finite type 
and $\fF _1 \in D(X\times Q_1)$ 
such that any object in $Q(G, \phi _0, \phi _1)$ is 
isomorphic to $\fF _{1,q}$ for some $q\in Q_1$. 
By the assumption and Proposition~\ref{purpose}, 
the generic flatness holds for $\pP ((\phi -1, \phi ])$
with $\phi \in \iI$. 
Thus the locus 
$$Q_1 ^{\circ}= \{ q\in Q_1 \mid \fF _q \in \pP ((\phi _0, \phi _1]) \},$$
is open because we have,
$$\pP ((\phi _0, \phi _1]) =\pP ((\phi _0, \phi _0 +1]) \cap 
\pP ((\phi _1 -1 , \phi _1]).$$
Now we can follow the same construction as in Proposition~\ref{quot}
and obtain 
$Q_2 =Q_1 ^{\circ}\times S$, $Q_3 \to Q_2$, $\widetilde{\eE} \in \Coh (Q_3)$,
and $\qQ (\eE, \phi _0, \phi _1) \subset \mathbb{V}(\widetilde{\eE})$
 as desired.
\end{proof}

The following is the main theorem in this section. 
\begin{thm}\label{condition}
Under the Assumption~\ref{nbd}, assume 
that for any $\sigma =(Z,\pP) \in \vV$,
Problem~\ref{pro} for $\aA =\pP ((0,1])$ and
 Problem~\ref{further}
are true.
 Then 
for any $\sigma \in \Stab^{\ast}(X)$, $\alpha \in \nN(X)$ and 
$\phi \in \mathbb{R}$, the stack $\mM ^{(\alpha, \phi)}(\sigma)$
is an Artin stack of finite type over $\mathbb{C}$. 
\end{thm}
Note that by Lemma~\ref{artinstack} and Lemma~\ref{opp}, it suffices to 
check Claim~\ref{op2} and Problem~\ref{further}. 
Also note that by Proposition~\ref{purpose} and Lemma~\ref{op3} (ii), the 
result holds for any $\sigma \in \vV$.
We divide the proof into some steps. 
\begin{step}\label{hide2}
The result holds for an algebraic stability condition
$\sigma =(Z, \pP) \in \overline{\vV}$. 
\end{step}
\begin{proof}
First we show Claim~\ref{op2} holds. 
For a smooth quasi-projective variety $S$ and 
$\eE \in \mM (S)$, assume the locus $S^{\circ}$ defined 
by (\ref{locus}) is non-empty. 
   Note that $S^{\circ}$ is dense in $S$ by Lemma~\ref{op3} (i). 
  Since $\sigma \in \overline{\vV}$, there exists $\sigma '=(Z', \pP ')
   \in \vV$
  and $\phi _i \in (\phi -1/2, \phi +1/2) \cap \iI$ for $0\le i\le 5$ such that
  $$\pP  (\phi) \subset \pP ' ((\phi _0, \phi _1]) \subset 
  \pP  ((\phi _2, \phi _3]) \subset \pP ' ((\phi _4, \phi _5]) \subset
  \pP  \left(\left(\phi -\frac{1}{2}, \phi +\frac{1}{2}\right]\right).$$
  For $E\in \pP  ((\phi -1/2, \phi +1/2])$, we denote by $\phi (E) \in 
  (\phi -1/2, \phi +1/2]$ the phase with respect to the stability 
  function $Z$. 
  By the assumption for $\sigma ' \in \vV$ and 
  Proposition~\ref{purpose}, there is an open subset $S_1 \subset S$
  on which $\eE _s \in \pP ' ((\phi _0, \phi _1])$. 
  Now we have 
  \begin{align*}
 & \{ s\in S \mid \eE _s \mbox{ is not semistable in }\sigma \} \\
  &= \{ s\in S \mid \mbox{there is a strict epimorphism }
  \eE _s \twoheadrightarrow F \mbox{ in }
  \pP  ((\phi _2, \phi _3]) \mbox{ with }
  \phi (F)< \phi (\eE _s) \} \\
  & \subset  \{ s\in S  \mid \mbox{there is a 
  strict epimorphism }
  \eE _s \twoheadrightarrow F \mbox{ in }\pP ' ((\phi _4, \phi _5]) 
  \mbox{ with }
  \phi (F)< \phi (\eE _s) \}.
  \end{align*}
  On the other hand, assume there is a strict epimorphism 
  $\eE _s \twoheadrightarrow F$ in $\pP ' ((\phi _4, \phi _5])$
  with $\phi (F) < \phi (\eE _s)$. 
 Then it
  is a surjection in $\pP  ((\phi -1/2, \phi +1/2])$,
   and $\phi (F) < \phi (\eE _s)$ implies $\eE _s$ is not semistable 
   in $\sigma$. 
   Thus we obtain, 
   \begin{align*}
 & \{ s\in S \mid \eE _s \mbox{ is not semistable in }\sigma \} \\
  & = 
   \{ s\in S  \mid \mbox{there is a 
  strict epimorphism }
  \eE _s \twoheadrightarrow 
  F \mbox{ in }\pP ' ((\phi _4, \phi _5]) \mbox{ with }
  \phi (F)< \phi (\eE _s) \}.
  \end{align*}
  Let 
  $$\pi _{\phi _4, \phi _5} \colon 
  \qQ(\eE, \phi _4, \phi _5) \lr S$$
   be the $S$-scheme 
  constructed in Lemma~\ref{ad} applied for $\sigma '$. Let 
  $\eE _{\qQ(\eE, \phi _4, \phi _5)} \to \fF$
   be the universal epimorphism on $X\times \qQ(\eE, \phi _4, \phi _5)$
   and define 
   $\qQ^{\circ}(\eE, \phi _4, \phi _5)$ to be the locus 
   $$\qQ^{\circ}(\eE, \phi _4, \phi _5)\cneq 
   \{ q\in \qQ(\eE, \phi _4, \phi _5) \mid 
   \phi (\fF _q) <\phi (\eE _q) \}.$$
   Since $q \mapsto \phi (\fF _q)$ is locally constant on 
   $\qQ(\eE, \phi _4, \phi _5)$, 
   $\qQ^{\circ}(\eE, \phi _4, \phi _5)$ is a union 
   of connected components of $\qQ(\eE, \phi _4, \phi _5)$, 
   in particular it is of finite type over $S$. 
   Let $\pi _{\phi _4, \phi _5}^{\circ}$ be the restriction of 
   $\pi _{\phi _4, \phi _5}$ to 
   $\qQ^{\circ}(\eE, \phi _4, \phi _5)$. Then 
   for a point $s\in S$, $s\in S^{\circ}$ if and only 
   if $s \notin \im \pi _{\phi _4, \phi _5}^{\circ}$. This implies 
   $S \setminus \im \pi _{\phi _4, \phi _5}^{\circ}$ is dense, 
   thus there exists an open subset $U\subset 
   S \setminus \im \pi _{\phi _4, \phi _5}^{\circ}$.
   
   Next we check that
   $M^{(\alpha, \phi)}(\sigma)$
    is bounded. Take $E\in M^{(\alpha, \phi)}(\sigma)$ and 
   let $F_i \in \pP '((\phi _0, \phi _1])$ for $1\le i\le n(E)$
    be the semistable factors of $E$ in $\sigma '$. Because $\sigma '$ is 
    algebraic and $\phi _1 -\phi _0 <1$, the map $
    E \mapsto n(E)$ is bounded on $M^{(\alpha, \phi)}(\sigma)$
     and 
     $$\{ Z'(F_i) \in \mathbb{C} \mid 1\le i\le n(E), E\in 
     M^{(\alpha, \phi)}(\sigma)\},$$
      is a finite set. Since we assume 
     that Problem~\ref{further} is true for $\sigma '$, 
      the set of objects
    $$\{ F_i \mid 1\le i\le n(E), E\in M^{(\alpha, \phi)}(\sigma) \},$$
     is bounded.
    Thus $M^{(\alpha, \phi)}(\sigma)$ is also bounded by Lemma~\ref{useasy}.
    
   \end{proof}

   \begin{step}\label{hide3}
   The result holds for any algebraic stability 
   condition $\sigma \in \Stab ^{\ast}(X)$. 
   \end{step}
   \begin{proof}
   Note that $\Phi \in \Auteq D(X)$ induces a 1-isomorphism, 
   $$\mM \ni E \longmapsto \Phi(E)\in \mM.$$
   Also note that an action of 
    $g\in \widetilde{\GL}^{+}(2, \mathbb{R})$ does not change the 
    set of semistable objects. Thus
   we have 
   $$\mM ^{(\alpha, \phi)}(\sigma) = \mM ^{(\alpha, \phi ')}(g(\sigma)),$$
   for some $\phi ' \in \mathbb{R}$. Hence if the result holds for 
   $\sigma \in \Stab ^{\ast}(X)$,
    then it also holds for $g\circ \Phi(\sigma)$ for 
   any $\Phi \in \Auteq D(X)$ and $g\in \widetilde{\GL}^{+}(2, \mathbb{R})$. 
   Thus the result holds for any algebraic stability 
   condition $\sigma \in \Stab ^{\ast}(X)$
    by Assumption~\ref{nbd} and Step~\ref{hide2}.     
  \end{proof}

   \begin{step}
    The result holds for any 
   $\sigma =(Z, \pP) \in \Stab ^{\ast}(X)$. 
   \end{step}
\begin{proof}
Let $\sigma \in \mathfrak{B}^{\circ}
 \subset \Stab ^{\ast}(X)$ be an open neighborhood 
of $\sigma$ such that its closure $\mathfrak{B}$ is compact. 
Let $\sS \subset D(X)$ be
$$\sS \cneq \{ E\in D(X) \mid 
E \mbox{ is of numerical type }\alpha \mbox{ and semistable in some }
\sigma ' \in \mathfrak{B} \}.$$
Then $\sS$ has bounded mass, 
hence by Assumption~\ref{nbd} and
Proposition~\ref{massfin}, there exists a finite 
number of codimension one walls $\{ \wW _{\gamma} \} _{\gamma \in \Gamma}$
which gives a wall and chamber structure on $\mathfrak{B}$. 
Let $\Gamma '\subset \Gamma$ be the subset which satisfies,
\begin{align}\label{ood}
\sigma \in \bigcap _{\gamma \in \Gamma '} \wW _{\gamma}
\setminus \bigcup_{\gamma \notin \Gamma '} \wW _{\gamma}.\end{align}
Let $\cC$ be the connected component of the right hand side 
of $(\ref{ood})$ which contains $\sigma$. 
Then if $E$ is of numerical type $\alpha$ and
semistable in $\sigma$, then 
it is semistable for any $\sigma '\in \cC$.  
We can take $\sigma '=(Z', \pP ')$ to be algebraic by Lemma~\ref{put}. 
Thus $\mM ^{(\alpha, \phi)}(\sigma) = \mM ^{(\alpha, \phi ')}(\sigma ')$
for some $\phi '$, and the result follows from Step~\ref{hide3}.

\end{proof}

\begin{rmk}\label{remark}\emph{
Note that Assumption~\ref{nbd} and Proposition~\ref{purpose}
also imply the following:
the set of $\sigma =(Z, \pP) \in \Stab ^{\ast}(X)$ such that 
$\pP ((\phi -1, \phi])$ satisfies the generic flatness for any 
$\phi \in \iI$ is dense in $\Stab ^{\ast}(X)$. }
\end{rmk}

\section{Semistable objects on K3 surfaces}\label{K3}
In this section we assume $X$ is a K3 surface or an 
abelian surface. 
The aim of this section is to show that the assumption in 
Theorem~\ref{condition} is satisfied in this case. 
\subsection{Mukai lattices and Mukai vectors}
Let $\NS ^{\ast}(X)$ be the Mukai lattice, 
$$\NS ^{\ast}(X) = \mathbb{Z} \oplus \NS (X) \oplus \mathbb{Z}.$$
For $v_i =(r_i, l_i,s_i)$ with $i=1,2$, its bilinear pairing is 
given by 
\begin{align}\label{pairing}
(v_1, v_2)=l_1\cdot l_2 -r_1s_2 -r_2 s_1. \end{align}
For an object $E\in D(X)$ its Mukai vector is defined as follows. 
\begin{align*}v(E) &=\ch (E)\sqrt{\td _X}\\
&= (r(E), c_1 (E), \ch _2 (E)+\epsilon \cdot r(E)).
\end{align*} 
Here $\epsilon =1$ if $X$ is a K3 surface and $\epsilon =0$
if $X$ is an abelian surface. Sending an object to its Mukai 
vector
gives an isomorphism, 
\begin{align}\label{muk}
v\colon \nN (X) \stackrel{\cong}{\lr} \NS ^{\ast}(X).\end{align}
Under the identification (\ref{muk}), the bilinear pairing 
$-\chi (E_1,E_2)$ on the left hand side goes to 
the pairing (\ref{pairing}). 
\subsection{Twisted Gieseker-stability and $\mu$-stability}
We recall the notion of 
twisted Gieseker-stability and $\mu$-stability
on the category of coherent sheaves $\Coh (X)$.
For the detail, one can consult~\cite{Hu}, \cite{MW}.  
Take $\lL, \mM \in \Pic(X)$, and suppose $\lL$ is ample. 
For $E\in \Coh (X)$
 one can 
write the twisted Hilbert polynomial as follows, 
$$
\chi (E\otimes\mM ^{-1} \otimes \lL ^n)= \sum _{i=0}^d a_i n^i,$$
for $a_i \in \mathbb{Q}$ and $a_d \neq 0$.  
 For $\omega= c_1 (\lL)$ and $\beta =c_1 (\mM)$, define the 
\textit{twisted reduced Hilbert polynomial} $P(E, \beta, \omega, n)$ 
to be
\begin{align}\label{twist}
P(E, \beta, \omega, n) &= 
\frac{\chi (E\otimes \mM ^{-1} \otimes \lL ^{n})}{a_d}\end{align}
When $\beta =0$, we simply 
write it $P(E, \omega, n)$.
Note that (\ref{twist}) is calculated by chern characters of 
$E$, $\mM$ and $\lL$. Thus by formally replacing the chern characters by 
their fractional, we can define $P(E, \beta, \omega, n)$ for 
$\mathbb{Q}$-divisors $\beta$ and $\omega$, and $E\in \nN (X)$.  
Explicitly when $v(E)=(r,l,s)$ with $r>0$, we have
\begin{align}\label{expli} (\ref{twist})= n^2 +
 \frac{2(l-r\beta)\cdot \omega}{r\omega ^2}n -
\frac{\{ l^2 -2rs -(l-r\beta)^2 \}}{r^2 \omega ^2}+ 
\frac{2\epsilon}{\omega ^2},
\end{align}
and $(\ref{twist})= n+ (s-\beta \cdot l)/\omega \cdot l$ when $r=0$, $l\neq 0$,
and $(\ref{twist})=s$ when $r=l=0$. 
Also for a torsion free sheaf $E$, define 
 $\mu _{\omega}(E) \in \mathbb{Q}$ to be
$$\mu _{\omega}(E)=\frac{l \cdot \omega}{r}.$$
\begin{defi}\emph{
For a pure sheaf $E\in \Coh (X)$, we say $E$ is $(\beta, \omega)$-twisted
(semi)stable if for any subsheaf $F\subsetneq E$ one has 
$$P(F, \beta, \omega, n) < P(E, \beta, \omega, n), \quad ( \mbox{resp} \le ),$$
for $n\gg 0$.
If $\beta =0$, we say simply $\omega$-Gieseker (semi)stable. 
Also a torsion free sheaf
 $E$ is $\mu _{\omega}$-(semi)stable if for any subsheaf 
$F\subsetneq E$ one has 
$$\mu _{\omega}(F) < \mu _{\omega}(E), \quad ( \mbox{resp} \le ).$$}
\end{defi}
There are notions of Harder-Narasimhan filtrations in both 
stability conditions~\cite{MW}.

\subsection{Stability conditions on K3 surfaces}
Here we recall the constructions of stability conditions on 
a K3 surface or an abelian surface $X$ studied in~\cite{Brs2}. 
Let $\beta, \omega$ be $\mathbb{Q}$-divisors on $X$ with 
$\omega$ ample. 
 For a torsion free sheaf $E \in \Coh (X)$, 
one has the Harder-Narasimhan filtration 
$$0= E_0 \subset E_1 \subset \cdots \subset E_{n-1}\subset E_n =E,$$
such that $F_i =E_i /E_{i+1}$ is $\mu _{\omega}$-semistable and 
$\mu _{\omega}(F_{i})> \mu _{\omega}(F_{i+1})$. 
Then define $\tT _{(\beta, \omega)} \subset \Coh (X)$ to be
the subcategory consists of sheaves whose torsion free parts 
have $\mu _{\omega}$-semistable Harder-Narasimhan factors 
of slope $\mu _{\omega}(F_i) > \beta \cdot \omega$. Also 
define $\fF _{(\beta, \omega)}
\subset \Coh (X)$ to be the subcategory consists of 
torsion free sheaves whose $\mu _{\omega}$-semistable factors have slope 
$\mu _{\omega}(F_i) \le \beta \cdot \omega$.
\begin{defi} \emph{
We define $\aA _{(\beta, \omega)}$ to be 
$$\aA _{(\beta, \omega)}=
\{ E\in D(X) \mid H^{-1}(E) \in \fF _{(\beta, \omega)}, H^0(E) \in 
\tT _{(\beta, \omega)}\}. $$} \end{defi}
\begin{rmk}\emph{
Note that different choices of $\beta$, $\omega$ may
define the same category $\aA _{(\beta, \omega)}$. 
For instance, we have $\aA _{(\beta, k\omega)}=\aA _{(\beta, \omega)}$ 
for $k\in \mathbb{Q}_{\ge 1}$. } \end{rmk}

We define $Z_{(\beta, \omega)}\colon \nN (X)\to \mathbb{C}$
by the formula.  
\begin{align}\label{formula}
Z_{(\beta, \omega)}(E) =( \exp (\beta +i\omega), v(E)).\end{align}
Explicitly if $v(E)=(r,l,s)$ and $r \neq 0$,  then (\ref{formula})
is written as 
\begin{align}\label{explicit}
Z_{(\beta, \omega)}(E)= \frac{1}{2r}\left( 
(l^2 -2rs) +r^2 \omega ^2 -(l-r\beta)^2 \right) +i(l-r\beta)\cdot \omega.
\end{align}
If $r=0$, (\ref{formula}) is written as $Z(E)=(-s+l\cdot \beta) +i(l\cdot 
\omega)$. 
We define $\sigma _{(\beta, \omega)}$ to be the pair 
$(Z_{(\beta, \omega)}, \aA _{(\beta, \omega)})$. 
\begin{prop}\emph{\bf{\cite[Lemma 6.2, Proposition 7.1]{Brs2}}}\label{const} 
The subcategory $\aA _{(\beta, \omega)}\subset D(X)$ is a heart of a 
bounded t-structure, and the pair $\sigma _{(\beta, \omega)}$
gives a stability condition on $D(X)$ if and only if for 
any spherical sheaf $E$ on $X$, one has $Z_{(\beta, \omega)}(E)\notin 
\mathbb{R}_{\le 0}$. This holds whenever $\omega ^2 >2$. 
\end{prop}
Let $\Stab ^{\ast}(X)$ be the connected component of 
$\Stab (X)$ which contains $\sigma _{(\beta, \omega)}$, and
define $\vV \subset \Stab ^{\ast}(X)$ to be
$$\vV =\{ \sigma _{(\beta, \omega)}\in \Stab ^{\ast}(X) \mid
\sigma _{(\beta, \omega)}
\mbox{ satisfies the assumption in Proposition~\ref{const}}.\}
$$
The following is stated in~\cite[Section 13]{Brs2}. 
\begin{thm}\label{stated}
The connected component $\Stab ^{\ast}(X)$ and the
subset $\vV \subset \Stab ^{\ast}(X)$ satisfy 
Assumption~\ref{nbd}.
\end{thm}
Finally we give the following useful lemma. 
\begin{lem}\label{useful}
(i) If $E\in D(X)$ satisfies $\Hom (E,E)=\mathbb{C}$, then 
$v(E)^2 \ge -2$. 

(ii) For $\mathbb{Q}$-divisors $\beta, \omega$ with $\omega$ ample and 
$m\in \mathbb{R}_{>0}$, 
the set of Mukai vectors 
$$\{ v\in \NS ^{\ast}(X) \mid v^2 \ge -2, \lvert ( \exp (\beta +i\omega), v)
\rvert \le m\},$$
is finite.  

(iii) For $E\in \nN (X)$, assume $v(E)=(0,l,s) \in \NS ^{\ast}(X)$
with $l\neq 0$. 
Then 
\begin{align}\label{onedim}
P(E, \beta, \omega, n)= n -\frac{\Ree Z_{(\beta, \omega)}(E)}{\Imm Z_{(\beta, \omega)}(E)} \in \mathbb{Q}[n].\end{align}
\end{lem}

(iv) For $E, E' \in \nN (X)$, $P(E, \beta, \omega, n)=P(E', \beta, \omega, n)$
if and only if 
$$\Imm \frac{Z_{(\beta, k\omega)}(E')}{Z_{(\beta, k\omega)}(E)}=0,$$
for infinitely many $k\in \mathbb{Q}$. 
\begin{proof}
(i) is proved in~\cite[Lemma 5.1]{Brs2}
 and (ii) is proved in~\cite[Lemma 8.2]{Brs2}.
 (iii) and (iv) follow easily from (\ref{expli}) and (\ref{explicit}).  
\end{proof}

\subsection{Generic flatness for $\aA _{(\beta, \omega)}$. }
Here we show the generic flatness in a special case. 

\begin{lem}\label{special}
Problem~\ref{pro} is true for $\aA = \aA _{(\beta, \omega)}$. 
\end{lem}
\begin{proof}
Let $S$ be a smooth projective variety over $\mathbb{C}$, 
$\lL \in \Pic (S)$ be an ample line bundle. Let us take 
$\eE\in \aA _S$. By the definition of $\aA _S$, we have 
$$\dR p_{\ast}(\eE\otimes \lL ^n ) \in \aA _{(\beta, \omega)},$$
for $n\gg 0$. In particular $\dR p_{\ast}(\eE\otimes \lL ^n )$
is concentrated in degree $[-1,0]$. Note that the following spectral 
sequence 
$$E_2 ^{i,j} =R^i p_{\ast}(H^j(\eE) \otimes \lL ^n ) \Rightarrow 
\dR ^{i+j}p_{\ast}(\eE\otimes \lL ^n),$$
degenerates for $n\gg 0$. Therefore $H^j (\eE)=0$ unless $j=-1$ or 
$0$. By~\cite[Theorem 2.3.2]{Hu}, there is an open subset 
$U\subset S$ and filtrations of coherent sheaves, 
$$
H^{-1}(\eE)_U =F^0 \supset F^1 \supset \cdots \supset F^k, \quad
H^0 (\eE)_U =T^0 \supset T^1 \supset \cdots \supset T^l,
$$
such that 
\begin{itemize}
\item Each $F^i$ and $T^i$ are flat sheaves on $U$. 
\item For $s\in U$, the filtrations
$$
H^{-1}(\eE)_s =F^0 _s \supset F^1 _s \supset \cdots \supset F^k _s, \quad
H^0 (\eE)_s =T^0 _s \supset T^1 _s \supset \cdots \supset T^l _s,
$$
are Harder-Narasimhan filtrations in $\omega$-Gieseker stability.
\end{itemize}
Note that $\eE _s \in \aA _{(\beta, \omega)}$ is equivalent to 
\begin{align}\label{numm}
\mu _{\omega}(F_s ^k) \le \beta \cdot \omega, \quad 
\mu _{\omega}(T_s ^0 /T_s ^1) > \beta \cdot \omega \mbox{ or } H^0(\eE)_s
 \mbox{ is 
torsion },\end{align}
and such points are dense in $S$ by Proposition~\ref{dense}. 
For
 each $i$, $s, s'\in U$, the coherent sheaves $F_s ^i$, $T_s ^i$
are numerically equivalent to $F_{s'}^i$, $T_{s'}^i$ respectively. 
Therefore (\ref{numm}) holds for any $s\in U$.  This implies  
$\eE_s \in \aA _{(\beta, \omega)}$ for any $s\in U$. 
\end{proof}

\subsection{Boundedness of semistable objects}
Next we check the boundedness of 
$M^{\alpha}(\sigma _{(\beta, \omega)})$,
where $\sigma _{(\beta, \omega)}\in \vV$. 
Let us prepare some notation and lemmas. For 
$E\in \aA _{(\beta, \omega)}$, let 
$$ H^0 (E) _{\rm{tor}} \subset H^0(E),$$
be the maximal 
torsion subsheaf of $H^0(E)$, and set 
$$H^0 (E)_{\rm{fr}} =H^0 (E)/H^0(E)_{\rm{tor}}.$$
Let 
\begin{align*}
& T_1, \cdots, T_{a(E)}\in \Coh (X), \\
& F_1, \cdots, F_{d(E)}, F_{d(E)+1}, \cdots, F_{e(E)} \in \Coh (X),
\end{align*}
be $\mu _{\omega}$-stable
 factors of $H^0(E)_{\rm{fr}}$, $H^{-1}(E)$ respectively. 
 Also let
$$ T_{a(E)+1}, \cdots, T_{b(E)}, T_{b(E)+1}, 
\cdots, T_{c(E)} \in \Coh (X)$$
  be $(\beta, \omega)$-twisted stable
 factors of $H^0(E)_{\rm{tor}}$.
 For the numbering, we set as follows. 
\begin{align*}
& \dim T_{i} =2 \quad (1\le i \le a(E)), \\
& \dim T_i =1 \quad (a(E) < i \le 
b(E)),  \\
& \dim T_i =0 \quad  (b(E)<i \le c(E)),\\
& \Imm Z_{(\beta, \omega)}(F_i [1])>0  \quad (1\le i\le d(E)), \\
& \Imm Z_{(\beta, \omega)}(F_i [1])=0 \quad (d(E)< i\le e(E)). 
\end{align*} 
Also for $\alpha \in \nN (X)$, define 
the set of objects 
$M^{\alpha}(\beta, \omega)$
to be 
$$M^{\alpha}(\beta, \omega)=
\{ E\in \aA _{(\beta, \omega)} \mid  \Imm Z_{(\beta, \omega)}(E)
\le \Imm Z_{(\beta, \omega)}(\alpha) \}.$$
We prepare the following lemma. 
\begin{lem}\label{con1}
The maps on $M^{\alpha}(\beta, \omega)$, 
$$E \longmapsto b(E), \quad E \longmapsto d(E),$$
are bounded. Furthermore the sets
\begin{align}\label{kakumei}
& \{ \Imm Z_{(\beta, \omega)}(T_i) \in \mathbb{Q} \mid 
1\le i\le c(E), E\in M^{\alpha}(\beta, \omega) \},\\
& \{ \Imm Z_{(\beta, \omega)}(F_i[1]) \in \mathbb{Q} \mid
1\le i\le e(E), E \in M^{\alpha}(\beta, \omega) \},
\end{align}
are finite sets. 
\end{lem}
\begin{proof}
For $E\in M^{\alpha}(\beta, \omega)$,
we have the inequality, 
$$\Imm Z_{(\beta, \omega)}(\alpha)\ge 
\Imm Z_{(\beta, \omega)}(E) = 
\sum_{i=1}^{b(E)} \Imm Z_{(\beta, \omega)}(T_i) +\sum _{i=1}^{d(E)}
\Imm Z_{(\beta, \omega)}(F_i[1]).$$
Note that each term of the above sum is positive. 
Noting that $\beta$ and $\omega$ are rational, we 
can conclude the result.
\end{proof}

The next step is to bound the real parts of $Z_{(\beta, \omega)}(T_i)$ 
and $Z_{(\beta, \omega)}(F_i[1])$. 
For the later use, we also give the bound of real part of 
$Z_{(\beta, k\omega)}(\ast)$ for $k\ge \mathbb{Q}_{\ge 1}$. 
\begin{lem}\label{fina}
There exist constants $C$, $C'$, $N$, which depend only on 
$\alpha$, $\beta$ and $\omega$ such that 
\begin{align}\label{constant1}
 &\frac{1}{k}\Ree Z_{(\beta, k\omega)}(T_i) \ge
  \Ree Z_{(\beta, \omega)}(T_i) \ge C \quad (1\le i\le a(E)), \\
\label{constant2} & \frac{1}{k}\Ree Z_{(\beta, k\omega)}(F_i[1])
\le \Ree Z_{(\beta, \omega)}(F_i [1])\le C' \quad (1\le i\le e(E)),
\end{align}
for any $E\in M^{\alpha}(\beta, \omega)$. 
\end{lem}
\begin{proof}
We give the proof of (\ref{constant1}). The proof of (\ref{constant2})
is similar. 
Denote 
 $$v(T_i)=(r_i, l_i, s_i) \in \mathbb{Z} \oplus \NS (X) \oplus \mathbb{Z}.$$
Note that $r_i >0$ for $1\le i\le a(E)$, and 
$$\Imm Z_{(\beta, \omega)}(T_i)=(l_i -r_i \beta )\cdot \omega, $$
which is bounded by Lemma~\ref{con1}. Thus
the Hodge index theorem implies that there exists a constant $C '' >0$
which depends only on $\alpha$, $\beta$ and $\omega$ such
that 
$$(l_i -r_i \beta)^2 
\le C''.$$
By Lemma~\ref{useful} (i), we have 
$$v(T_i)^2 =l_i ^2 -2r_i s_i \ge -2.$$
Hence for $1\le i\le a(E)$, we have 
\begin{align*}
\Ree Z_{(\beta, \omega)}(T_i) &=
\frac{1}{2r_i}\left( (l_i ^2 -2r_i s_i )+r_i ^2 \omega ^2 
-(l_i -r_i \beta)^2 \right)  \\
& \ge \frac{1}{2}r_i \omega ^2 -\frac{2+C''}{2r_i} \\
& > -\frac{2+C''}{2}.
\end{align*}
Similarly we have 
\begin{align*}
\frac{1}{k}\Ree Z_{(\beta, k\omega)}(T_i) -\Ree Z_{(\beta, \omega)}(T_i)
\ge \frac{1}{2}(k-1)\omega ^2 +\left(\frac{1}{k}-1\right) \frac{2+C''}{2},
\end{align*}
Thus one can find a desired $N>0$. 

\end{proof}

Finally we give the following preparation. 
\begin{lem}\label{asu}
Let $\sS$ be a subset of $M^{\alpha}(\beta, \omega)$.

(i)
Assume 
$$E\longmapsto \Ree Z_{(\beta, \omega)}H^{0}(E)_{\rm{fr}},$$
is bounded above on $\sS$. Then the following set,
\begin{align}\label{bel}
\{ v(T_i) \in \NS ^{\ast}(X) \mid 1\le i\le a(E), E\in \sS \},
\end{align}
is a finite set. 

(ii)
Assume 
$$E\longmapsto \Ree Z_{(\beta, \omega)}H^{-1}(E),$$
is bounded below on $\sS$. Then the following set,
\begin{align}\label{abo}
\{ v(F_i) \in \NS ^{\ast}(X) \mid 1\le i\le e(E), E\in \sS \},
\end{align}
is a finite set. 
\end{lem}
\begin{proof}
We show (ii). The proof of (i) is similar and leave it to the 
reader. 
For $E\in \sS$, we have 
\begin{align}
\label{sen} & \Ree Z_{(\beta, \omega)}(H^{-1}(E)[1])=\sum _{i=1}^{d(E)}
\Ree Z_{(\beta, \omega)}(F_i[1])
+ \sum _{i=d(E)+1}^{e(E)}\Ree Z_{(\beta, \omega)}(F_i[1]).
\end{align}
Note that $Z_{(\beta, \omega)}(F_i[1])\in \mathbb{R}_{<0}$ for 
$d(E)<i\le e(E)$, and 
$\Ree Z_{(\beta, \omega)}(F_i[1])$ is bounded above
for $1\le i\le d(E)$
 by 
Lemma~\ref{fina}. Furthermore $E\mapsto d(E)$ is bounded by Lemma~\ref{con1}. 
Therefore  
the map $E \mapsto e(E)$ is bounded and 
the following set is a finite set: 
\begin{align*}
\{ \Ree Z_{(\beta, \omega)}(F_i[1]) \in \mathbb{Q} \mid 
1\le i\le e(E), E\in \sS \}.
\end{align*}
Then combined with Lemma~\ref{con1}, the
following set is a finite set:
\begin{align}\label{eien}
 \{ Z_{(\beta, \omega)}(F_i[1]) \in \mathbb{C} \mid 
1\le i\le e(E), E\in \sS \}.
\end{align}
By the finiteness 
of (\ref{eien}) and Lemma~\ref{useful} (i), (ii), 
the set (\ref{abo})
is also finite. 

\end{proof}

Now we can show the following. 
\begin{prop}\label{bousemi}
Problem~\ref{further} is 
true for any 
$\sigma _{(\beta, \omega)}\in \vV$. 
\end{prop}
It is enough to show the boundedness of
$$M^{\alpha}(\sigma _{(\beta, \omega)})
=\{ E\in \aA _{(\beta, \omega)}\mid E \mbox{ is of numerical type }\alpha
\mbox{ and semistable in }\sigma _{(\beta, \omega)} \}.$$
Note that we have $M^{\alpha}(\sigma _{(\beta, \omega)})\subset 
M^{\alpha}(\beta, \omega)$. 
Let $\tT$, $\tT '$ and $\fF$ be the sets of objects, 
\begin{align*}  \tT &=\{H^0 (E)_{\rm{fr}} \in \Coh (X) \mid E\in 
M^{\alpha}(\sigma _{(\beta, \omega)})\}, \\
 \tT ' &= \{ H^0 (E)_{\rm{tor}} \in \Coh (X) \mid
 E\in M^{\alpha}(\sigma _{(\beta, \omega)}) \},\\
 \fF &=\{ H^{-1}(E) \in \Coh (X) \mid E\in 
 M^{\alpha}(\sigma _{(\beta, \omega)})\}.
 \end{align*}
 By Lemma~\ref{useasy}, it
  suffices to show that each $\tT$, $\tT'$ and $\fF$ are bounded. 
 We divide the proof into two steps. 

\begin{sstep}
The sets of objects $\tT$, $\fF$ are bounded. 
\end{sstep}
\begin{proof}
Take $E\in M^{\alpha}(\sigma _{(\beta, \omega)})$. Note that
we have the exact sequence in $\aA _{(\beta, \omega)}$, 
$$0 \lr H^{-1}(E)[1] \lr E \lr H^0 (E) \lr 0,$$
and a surjection 
$H^0(E) \twoheadrightarrow H^0(E)_{\rm{fr}}$
in $\aA _{(\beta, \omega)}$. Thus we have 
$$\Imm Z_{(\beta, \omega)}(H^{-1}(E)[1])\le 
\Imm Z_{(\beta, \omega)}(\alpha), \quad
 \Imm Z_{(\beta, \omega)}(H^0 (E)_{\rm{fr}})
\le \Imm Z_{(\beta, \omega)}(\alpha),$$
and
the semistability of $E$ implies 
$$\phi (H^{-1}(E)[1]) \le \phi (E) \le \phi (H^0(E)_{\rm{fr}}).$$
Therefore if we consider the maps on 
$M^{\alpha}(\sigma _{(\beta, \omega)})$, 
\begin{align}\label{bengo}
& E\longmapsto \Ree Z_{(\beta, \omega)}(H^0(E)_{\rm{fr}})\in \mathbb{Q}, \\
\label{sumita} 
& E \longmapsto \Ree Z_{(\beta, \omega)}(H^{-1}(E)[1])\in \mathbb{Q},
\end{align}
then (\ref{bengo}) is bounded above and (\ref{sumita}) is 
bounded below. 
Thus one can apply Lemma~\ref{asu} and conclude that the sets
\begin{align*}
& \{ v(T_i) \in \NS ^{\ast}(X) \mid 
1\le i\le a(E), E\in M^{\alpha}(\sigma _{(\beta, \omega)})\},\\
&  \{ v(F_i[1]) \in \NS ^{\ast}(X) \mid 
1\le i\le e(E), E\in M^{\alpha}(\sigma _{(\beta, \omega)})\},
\end{align*}
are finite sets. Since the set of 
Gieseker-stable sheaves with a fixed Mukai vector is bounded, 
(see~\cite{Hu}) 
the sets of sheaves 
\begin{align*}
& \{ T_i \in \Coh (X) \mid 
1\le i\le a(E), E\in M^{\alpha}(\sigma _{(\beta, \omega)})\},\\
&  \{ F_i \in \Coh (X) \mid 
1\le i\le e(E), E\in M^{\alpha}(\sigma _{(\beta, \omega)})\},
\end{align*}
are bounded. Thus $\tT$ and $\fF$ are also bounded by Lemma~\ref{useasy}
\end{proof}
\begin{sstep}
The set of sheaves $\tT'$ is bounded. 
\end{sstep}
\begin{proof}
For $a(E)<i\le b(E)$
we may assume 
$P(T_{i},\beta, \omega,n)> P(T_{i+1},\beta, \omega,n)$.
 Hence by Lemma~\ref{useful} (iii) we have
 \begin{align}\label{num}\phi(T_{a(E)+1})> \cdots > \phi(T_{b(E)}).
 \end{align}
 Note that there is an exact sequence
 $$0 \lr T' \lr H^0(E)_{\rm{tor}} \lr T_{b(E)} \lr 0,$$
  both in $\Coh (X)$ and $\aA _{(\beta, \omega)}$. 
  Let $H^0(E)/T' \in \Coh (X)$ be the cokernel of the inclusion,
  $$T' \hookrightarrow H^0(E)_{\rm{tor}} 
\hookrightarrow H^0(E),$$
in $\Coh (X)$. 
Then the following composition,
$$E \twoheadrightarrow H^0(E) \to H^0(E)/T',$$
is a surjection in $\aA _{(\beta,\omega)}$.  
Thus 
we have 
$$\Imm Z_{(\beta, \omega)}(H^0(E)/T')\le \Imm Z_{(\beta, \omega)}(E),$$
and
the 
semistability of $E$ implies 
$\phi (E)\le \phi(H^0(E)/T')$. Hence the map
\begin{align*}
E\longmapsto & \Ree Z_{(\beta, \omega)}(H^0(E)/T'), \\
&= \Ree Z_{(\beta, \omega)}(T_{b(E)})+ \Ree Z_{(\beta, \omega)}(H^0(E)_{\rm{fr}}),
\end{align*}
is bounded above. Since $\tT$ is bounded, it follows that 
$\Ree Z_{(\beta, \omega)}(T_{b(E)})$ is also bounded above.
Hence by Lemma~\ref{con1} and (\ref{num}), there is a constant 
$C'''$ (which depends only on $\alpha$, $\beta$ and $\omega$) such that 
\begin{align}\label{hen}
\Ree Z_{(\beta, \omega)}(T_i) \le C''', \quad (a(E)<i\le b(E)).
\end{align} 
On the other hand we have 
\begin{align}\label{hey}
\Ree Z_{(\beta, \omega)}(H^0(E)_{\rm{tor}})=
\sum _{i=a(E)+1}^{b(E)}\Ree Z_{(\beta, \omega)}(T_i)+ 
\sum _{i=b(E)+1}^{c(E)} \Ree Z_{(\beta, \omega)}(T_i).
\end{align}
Note that $E\mapsto \Ree Z_{(\beta, \omega)}(H^0(E)_{\rm{tor}})$ is bounded 
on $M^{\alpha}(\sigma _{(\beta, \omega)})$
because $\tT$ and $\fF$ are bounded. Thus the boundedness of (\ref{hey})
together with (\ref{hen}) and $\Ree Z_{(\beta, \omega)}(T_i) 
\in \mathbb{R}_{<0}$ for $b(E)<i\le c(E)$ show that the set, 
\begin{align}\label{sou}
\{ \Ree Z_{(\beta, \omega)}(T_i) \in \mathbb{Q} \mid 
a(E)< i\le c(E), E\in M^{\alpha}(\sigma _{(\beta, \omega)})\},
\end{align}
is a finite set, and $E\mapsto c(E)$ is bounded.
 Hence by Lemma~\ref{con1}, the finiteness of (\ref{sou}), 
and using Lemma~\ref{useful} (ii), we conclude that
the set 
$$\{ v(T_i)\in \NS ^{\ast}(X) \mid a(E)< i\le c(E), E\in 
M^{\alpha}(\sigma _{(\beta, \omega)})\},$$
is a finite set. Again the set of sheaves, 
$$\{ T_i \in \Coh (X) \mid a(E)<i\le c(E), E\in 
M^{\alpha}(\sigma _{(\beta, \omega)}) \},$$
is bounded, thus $\tT'$ is also bounded by Lemma~\ref{useasy}.
\end{proof}

Combined with the result in the previous section, we obtain 
the following. 

\begin{thm}\label{sekai}
Let $X$ be a K3 surface or an abelian surface. Then
for any $\sigma \in \Stab ^{\ast}(X)$, $\alpha \in \nN (X)$, 
$\phi \in\mathbb{R}$, the stack $\mM ^{(\alpha, \phi)}(\sigma)$ 
is an Artin stack of finite type over $\mathbb{C}$. 
\end{thm}
\begin{proof}
This follows from Theorem~\ref{condition}, Theorem~\ref{stated},
 Lemma~\ref{special}, 
and Proposition~\ref{bousemi}. 
\end{proof}

\section{Invariants counting semistable objects}\label{counting}
In this section, $X$ is a K3 surface or an abelian surface,  
 $\mM$ is the moduli stack of objects $E\in D(X)$ with 
 $\Ext ^{<0}(E,E)=0$ as in the previous section. 
The aim in this section is to introduce and study the
 the invariants, as an analogue of the work~\cite{Joy4}. 
\subsection{Stack functions}
Let $\dD$ be an Artin stack over $\mathbb{C}$. 
Following D.~Joyce's work~\cite{Joy5}, we introduce the notion of
\textit{stack functions} on $\dD$. 
For the detail, one can consult~\cite[Section 3]{Joy5}.
Let us consider pairs $(\rR, \rho)$, where $\rR$ is
an Artin $\mathbb{C}$-stack of finite type over $\mathbb{C}$
with affine geometric stabilizers and $\rho \colon \rR \to \dD$
is a 1-morphism. 
We say two pairs $(\rR, \rho)$, $(\rR', \rho')$ \textit{equivalent}
if there exists a 1-isomorphism $\tau \colon \rR \to \rR '$ such that 
$\rho '\circ \tau$ is 2-isomorphic to $\rho$. 
\begin{defi}\emph{
Define the $\mathbb{Q}$-vector space 
$\SF (\dD)$ to be 
$$\SF (\dD) \cneq \bigoplus _{(\rR, \rho)}\mathbb{Q}[(\rR, \rho)]/\sim.$$
Here $[(\rR, \rho)]$ is an equivalence class of $(\rR, \rho)$ and 
the relation $\sim$ is generated by 
$$[(\rR, \rho)]=[(\rR ^{\dag}, \rho |_{\rR ^{\dag}})]+
[(\rR \setminus \rR ^{\dag}, 
\rho |_{\rR \setminus \rR ^{\dag}})],$$
where $\rR ^{\dag}$ is a closed substack of $\rR$. }
\end{defi}
For $\rho \colon \rR \to \dD$ and $\rho ' \colon \rR ' \to \dD$, 
there is a notion of fiber product~\cite[Definition 2.10]{Joy1}, 
$$\xymatrix{
\rR \times _{\rho, \dD, \rho '}\rR ' \ar[r]
\ar[d]_{\pi _{\rR '}}
& \rR \ar[d]^{\rho} \\
\rR ' \ar[r]^{\rho '} & \dD.}$$
As in~\cite[Definition 3.1]{Joy5}, we can define a $\mathbb{Q}$-bilinear 
product
$\SF (\dD) \times \SF (\dD) \to\SF (\dD)$ 
 by the formula, 
$$[(\rR, \rho)]\cdot [(\rR ', \rho ')]=
[(\rR \times _{\rho, \dD, \rho '} \rR ', \rho ' \circ \pi _{\rR '})].$$
Let $\Pi \colon \dD \to \cC$ be a 1-morphism of Artin $\mathbb{C}$-stacks.
Then define the \textit{push-forward}
 $\Pi _{\ast}\colon \SF (\dD) \to \SF (\cC)$
by 
$$\Pi _{\ast}\colon \sum _{i=1}^{m}c_i [(\rR _i, \rho _i)] 
\longmapsto \sum _{i=1}^{m}c_i [(\rR _i, \Pi \circ \rho _i)].$$
If $\Pi$ is of finite type, one can define the 
\textit{pull-back} $\Pi ^{\ast}\colon \SF (\cC) \to \SF (\dD)$,
$$\Pi ^{\ast}\colon \sum _{i=1}^{m}c_i [(\rR _i, \rho _i)] 
\longmapsto \sum _{i=1}^{m}c_i [(\rR _i \times _{\rho _i,\cC,\phi}\dD,
 \pi _{\dD})].$$ 
 The \textit{tensor product} $\otimes \colon \SF (\dD) \times \SF (\cC) 
 \to \SF (\dD \times \cC)$ is 
 $$ \left( \sum _{i=1}^{m}c_i [(\rR _i, \rho _i)]\right) \otimes 
  \left( \sum _{i=1}^{m'}c_i '[(\rR _i ', \rho _i ')] \right) 
  =\sum _{i,j}c_i d_j  [(\rR _i \times \rR _i ', \rho _i \times \rho _i ')]).
  $$
One can consult~\cite[Definition 3.1]{Joy5} for the detail of 
these definitions. 
  For a substack $i\colon \dD ^{\circ}\hookrightarrow \dD$, 
  we write $[(\dD, i)]$ as $[\dD ^{\circ}\hookrightarrow \dD]$. 
If $X$ is a K3 surface or an abelian surface, we have shown in 
Theorem~\ref{sekai} that
the stack $\mM ^{(\alpha, \phi)}(\sigma)$ is an open substack 
of $\mM$ and it is of finite type. 
\begin{defi}\emph{
For $\sigma \in \Stab ^{\ast}(X)$, $\alpha \in \nN (X)$ and 
$\phi \in \mathbb{R}$, 
we define $\delta ^{(\alpha, \phi)}(\sigma)$ to be 
$$\delta ^{(\alpha, \phi)}(\sigma) =[\mM ^{(\alpha, \phi)}(\sigma)
\hookrightarrow \mM]
\in \SF (\mM).$$}
\end{defi}
 
 \subsection{Ringel-Hall algebras}\label{Ringel}
Take an algebraic stability condition $\sigma =(Z, \pP)\in \Stab (X)$
and let $\aA _{\phi} =\pP ((\phi -1,\phi])$ for $\phi \in \iI$.  
Assume the generic flatness holds for $\aA _{\phi}$. Then 
the stack of objects in $\aA _{\phi}$ is an open substack of $\mM$, thus 
in particular it is an Artin stack over $\mathbb{C}$. 
We denote it by $\mathfrak{Obj}\aA _{\phi} \subset \mM$.
Following~\cite[Definition 5.1]{Joy2}, 
we introduce the associative multiplication 
$\ast$ on $\SF (\mathfrak{Obj}\aA _{\phi})$ based on \textit{Ringel-Hall algebras}. 
Let $\mathfrak{Obj}(\aA _{\phi}, n)$ be the moduli stack of filtrations, 
\begin{align}
\label{saki}
0=E_0 \subset E_1 \subset E_2 \subset \cdots \subset E_n,\end{align}
with $E\in \aA _{\phi}$. It is shown in~\cite[Theorem 8.2]{Joy1} that 
$\mathfrak{Obj}(\aA _{\phi}, n)$ is an Artin stack of locally 
finite type over $\mathbb{C}$. 
We have the following 1-morphisms, 
$$\prod _{i=1} ^n \mathfrak{Obj}\aA _{\phi} \stackrel{\prod _{i=1}^n p_i}
{\longleftarrow}
\mathfrak{Obj}(\aA _{\phi}, n) \stackrel{\Pi _n}{\longrightarrow}
 \mathfrak{Obj}\aA _{\phi}.$$
 Here $p_i \colon \mathfrak{Obj}(\aA _{\phi}, n) \to \mathfrak{Obj}\aA _{\phi}$ is 
 defined to be 
 $$(0=E_0 \subset E_1 \subset E_2 \subset \cdots \subset E_n)
 \longmapsto F_i \cneq E_i /E_{i-1},$$
 and $\Pi _n \colon \mathfrak{Obj}(\aA _{\phi}, n)
  \to \mathfrak{Obj}\aA _{\phi}$
 is defined to be 
 $$(0=E_0 \subset E_1 \subset E_2 \subset \cdots \subset E_n)
 \longmapsto E_n.$$
 It is shown in~\cite[Theorem 8.4]{Joy1} that $\prod _{i=1}^n p_i$ is of 
 finite type, thus one can define its pull-back. 
 One has the following diagram, 
 $$\xymatrix{
 \SF(\mathfrak{Obj}\aA _{\phi}) \times \SF (\mathfrak{Obj}\aA _{\phi}) 
 \ar[dr]
 \ar[d]_{\otimes} & & \\
 \SF(\mathfrak{Obj}\aA _{\phi} \times \mathfrak{Obj}\aA _{\phi}) 
 \ar[r]^{(p_1 \times p_2)^{\ast}} & 
 \SF(\mathfrak{Obj}(\aA _{\phi},2)) \ar[r]^{\Pi _{2\ast}} & \SF(\mathfrak{Obj}\aA _{\phi}),}$$
 \begin{defi}\emph{
We define a bilinear operation $\ast \colon \SF (\mathfrak{Obj}\aA _{\phi}) \times 
\SF (\mathfrak{Obj}\aA _{\phi}) \to \SF (\mathfrak{Obj}\aA _{\phi})$ to be 
$$f\ast g= \Pi _{2\ast}((p_1 \times p_2 )^{\ast}(f\otimes g)).$$}
\end{defi}
It is shown in~\cite[Theorem 5.2]{Joy2} that $\ast$ is associative and 
$\SF (\mathfrak{Obj}\aA _{\phi})$ is a $\mathbb{Q}$-algebra with identity
$\delta _{[0]} \cneq [(0\hookrightarrow \mathfrak{Obj}\aA _{\phi})]$.
In fact we have 
\begin{align}
& (f\ast g)\ast h =f \ast (g\ast h) =\Pi _{3\ast}
((p_1 \times p_2 \times p_3)^{\ast}(f\otimes g\otimes h)),\\
\label{asso}
& f_{1}\ast \cdots \ast f_n = \Pi _{n\ast}(\prod _{i=1}^n p_i)^{\ast}
(f_1 \otimes \cdots \otimes f_n ).\end{align}
One can consult~\cite[Section 5]{Joy2} for the detail 
of the algebra $(\SF (\mathfrak{Obj}\aA _{\phi}), \ast)$. 
For an interval $ I \subset \mathbb{R}$, set $C^{\sigma}(I)\subset
\nN (X)$ to be 
$$C^{\sigma}(I)= \im (\pP (I) \to \nN (X)) \setminus \{0\}
\subset \nN (X).$$
\begin{defi}\emph{
For $\alpha _1, \cdots, \alpha _n \in
 C^{\sigma}((\phi -1,\phi])$, 
we define the substack 
$$\mM (\{ \alpha _i \}_{1\le i\le n}, \aA _{\phi}, \sigma) \subset 
\mathfrak{Obj}(\aA _{\phi},n)$$ to be the stack of filtrations (\ref{saki})
such that $F_i =E_i /E_{i-1}$ is semistable in $\sigma$ and of numerical type 
$\alpha _i$. }
\end{defi}
Note that for
$\alpha \in C^{\sigma}((\phi -1,\phi])$, there is a 
unique phase $\phi (\alpha)\in (\phi -1, \phi]$
 with respect to the stability 
 function $Z$. 
 Also note that the
  element $\delta ^{(\alpha, \phi(\alpha))}(\sigma)\in \SF (\mM)$
 is regarded as the element of $\SF (\mathfrak{Obj}\aA _{\phi})$. 
 \begin{lem}
 For $\alpha _1, \cdots, \alpha _n \in C^{\sigma}((\phi -1, \phi])$, 
 we have the following equality 
 in $\SF (\mathfrak{Obj}\aA _{\phi})$,
  \begin{align}\label{kowareta}\Pi _{n\ast}
[(\mM (\{ \alpha _i\}_{1\le i\le n}, \aA _{\phi}, \sigma) \hookrightarrow 
\mathfrak{Obj}(\aA _{\phi}, n)]=
\delta ^{(\alpha _1, \phi (\alpha _1))}
(\sigma)\ast \cdots \ast \delta ^{(\alpha _n, 
\phi (\alpha _n))}(\sigma).\end{align}
\end{lem}
\begin{proof}
By the definition we have  
\begin{align}\label{mieru}
[\mM (\{ \alpha _i\}_{1\le i\le n}, \aA _{\phi}, \sigma) \hookrightarrow 
\mathfrak{Obj}(\aA _{\phi}, n)]=
(\prod _{i=1}^n p_i)^{\ast}[\prod \mM ^{(\alpha _i, \phi (\alpha _i))}(\sigma)
\hookrightarrow \prod \mathfrak{Obj}\aA _{\phi}],\end{align}
in $\SF (\mathfrak{Obj}(\aA _{\phi},n))$. 
Thus it is enough to apply $\Pi _{n\ast}$ to (\ref{mieru}) and
use (\ref{asso}). 
\end{proof}

\subsection{Motivic invariants of Artin stacks}
Let $K(\Var)$ be the Grothendieck ring of quasi-projective 
varieties. This is a $\mathbb{Z}$-module generated by the isomorphism 
classes of quasi-projective varieties $[X]$, and relations 
$[X]=[Y]+[X\setminus Y]$ for closed subschemes $Y\subset X$. 
The formula $[X]\cdot [X']=[X\times X']$ extends to a ring 
structure on $K(\Var)$. Suppose $\Lambda$ is a commutative
$\mathbb{Q}$-algebra and $\Upsilon$ is a ring homomorphism, 
\begin{align}\label{mot}
\Upsilon \colon K(\Var) \lr \Lambda. \end{align}
Write $l= \Upsilon (\mathbb{A}^1) \in \Lambda$. We assume $l$ and 
$l^k -1$ are invertible in $\Lambda$ for $k\ge 1$. This assumption 
is required for the value 
$$\Upsilon (\GL(m, \mathbb{C}))=l^{m(m-1)/2}\prod _{k=1}^m (l^k-1),$$
to be invertible in $\Lambda$.
\begin{exam}\label{virtual}\emph{
 We can take $\Lambda =\mathbb{Q}(z)$
and $\Upsilon([X])=P(X;z)$ the \textit{virtual Poincare polynomial}
of $X$. When $X$ is smooth and projective, $P(X;z)$ is the usual 
Poincare polynomial $\sum _{k=0}^{\dim X}b^k (X)z^k$. }
\end{exam}
An algebraic $\mathbb{C}$-group $G$ is called \textit{special}
if every principal $G$-bundle is locally trivial. It is shown 
in~\cite[Lemma 4.6]{Joy5} that if $G$ is special then 
$\Upsilon ([G])$ is invertible in $\Lambda$. 
\begin{thm}\emph{\bf{\cite[Theorem 4.9]{Joy5}}} 
Under the above situation, there exists a unique morphism 
of $\mathbb{Q}$-algebras, 
$$\Upsilon ' \colon \SF (\Spec \mathbb{C}) \lr \Lambda,$$
such that if $G$ is a special algebraic $\mathbb{C}$-group
which acts on a quasi-projective variety $X$, then 
$\Upsilon '([X/G])=\Upsilon ([X])/\Upsilon ([G])$. 
\end{thm}
Let $\Pi \colon \mM \to \Spec \mathbb{C}$ be the structure 
morphism. Given a motivic invariant $\Upsilon \colon K(\Var) \to \Lambda$
as in (\ref{mot}), 
we have the following maps, 
\begin{align}\label{signal}
\Upsilon ' \circ \Pi _{\ast} \colon \SF (\mM) \lr \SF (\Spec \mathbb{C}) 
\lr \Lambda.\end{align}

\begin{defi}\label{invariant}\emph{
Take $\sigma =(Z, \pP)\in \Stab ^{\ast}(X)$ and $\alpha \in \nN (X)$. 
We define $I^{\alpha}(\sigma)\in \Lambda$ as follows. If 
$Z(\alpha)=0$, we set $I^{\alpha}(\sigma)=0$. Otherwise
take $\phi \in \mathbb{R}$ which satisfies
$Z(\alpha)\in \mathbb{R}_{>0}e^{i\pi \phi}$, and define
$I^{\alpha}(\sigma)$ to be 
$$I^{\alpha}(\sigma)= \Upsilon ' \circ \Pi _{\ast}\delta ^{(\alpha, \phi)}
(\sigma) \in \Lambda.$$}
\end{defi}
The definition of $I^{\alpha}(\sigma)$ is an
 analogue of~\cite[Definition 6.1]{Joy4}.  
It is clear that the definition of $I^{\alpha}(\sigma)$ does not 
depend on a choice of $\phi$. Then as an analogue
 of~\cite[Definition 6.22]{Joy4}, 
we introduce the invariant $J^{\alpha}(\sigma) \in \Lambda$. 
\begin{defi}\emph{
We define $J^{\alpha}(\sigma)\in \Lambda$ as follows. 
If $Z(\alpha)=0$, we set $J^{\alpha}(\sigma)=0$. Otherwise 
choose $\phi \in \mathbb{R}$ which satisfies $Z(\alpha)\in 
\mathbb{R}_{>0}e^{i\pi \phi}$, and define $J^{\alpha}(\sigma)$ to be
\begin{align}\label{yomi}J^{\alpha}(\sigma)\cneq 
\sum _{
\alpha _1 + \cdots + \alpha _n =\alpha }
l^{-\sum _{j>i}\chi (\alpha _j, \alpha _i)}\frac{(-1)^{n-1}(l-1)}{n}
\prod _{i=1} ^n I^{\alpha _i}(\sigma) \in \Lambda,\end{align}
where $\alpha _i \in C^{\sigma}(\phi)$ for all $1\le i\le n$.  }
\end{defi}
Again the definition of $J^{\alpha}(\sigma)$ does not 
depend on a choice of $\phi$.
\begin{rmk}\label{stabilizers}\emph{
Suppose that any $E\in M^{(\alpha, \phi)}(\sigma)$ is 
stable. This occurs whenever $\alpha \in \nN (X)$ is primitive 
and $\sigma$ is not contained in a wall in Proposition~\ref{massfin}. 
In this case, 
$J^{\alpha}(\sigma)$ coincides with 
$(l-1)I^{\alpha}(\sigma)$. 
Furthermore any $E\in M^{(\alpha, \phi)}(\sigma)$
satisfies $\Hom (E, E)=\mathbb{C}$. Hence by Inaba's work~\cite{Inaba} 
and the openness of stability proved in Section~\ref{K3}, 
there is an algebraic space $\widetilde{M}^{(\alpha, \phi)}(\sigma)$
which parameterizes the objects in $M^{(\alpha, \phi)}(\sigma)$. 
Hence in this case we have 
$$J^{\alpha}(\sigma)=\Upsilon (\widetilde{M}^{(\alpha, \phi)}(\sigma))
\in \Lambda.$$
The factor $l-1 =\Upsilon (\mathbb{C}^{\ast})$ is required to cancel 
out the contributions of the stabilizers $\Aut (E)\cong \mathbb{C}^{\ast}$. }
\end{rmk}
\begin{rmk}\emph{
Suppose $\Lambda =\mathbb{Q}(z)$
 and $\Upsilon$ be as in Example~\ref{virtual}. 
Under the assumption in Remark~\ref{stabilizers}, 
we can define the invariant
 $J^{\alpha}(\sigma)|_{z=-1}\in \mathbb{Q}$, as a \textit{virtual euler 
 number} of the moduli space. However in general, 
 we don't know whether the denominator of $J^{\alpha}(\sigma)\in \mathbb{Q}(z)$
 is divided by $z+1$ or not. So at this time, we do not define the 
 invariant in $\mathbb{Q}$ in this way. 
 Also see Remark~\ref{motivation} below. }
\end{rmk}
We have to check the following. 
\begin{lem}
The sum (\ref{yomi}) is a finite sum.
\end{lem}
\begin{proof}
 Note that for any $\varepsilon >0$ 
there is an algebraic stability condition
$\sigma '=(Z',\pP ')\in \Stab ^{\ast}(X)$ such that
$\pP (\phi)\subset \pP '((\phi -\varepsilon, \phi +\varepsilon))$. 
Thus the possibilities of $n$ in the sum 
(\ref{yomi}) is finite. 
 Let us assume $\prod I^{\alpha _i}(\sigma) \neq 0$ in the sum (\ref{yomi}).
 Then there are objects $E_i \in \pP (\phi)$ of numerical type $\alpha _i$. 
 By taking stable factors of $E_i$ and
using Lemma~\ref{useful}, one can check that
  the possibilities for $\alpha _i$ is also finite. 
\end{proof}

\subsection{The algebra $A(\aA _{\phi}, \Lambda, \chi)$}
Let $\sigma =(Z,\pP)\in \Stab ^{\ast}(X)$ and 
$\aA _{\phi} =\pP ((\phi -1,\phi])$ be as in (\ref{Ringel}).  
We introduce the $\Lambda$-algebra $A(\aA _{\phi}, \Lambda, \chi)$. 
For the detail, see~\cite[Section 6]{Joy2}. 
\begin{defi}\emph{\bf{\cite[Definition 6.3]{Joy2}}} \emph{
We define the $\Lambda$-algebra $A(\aA _{\phi}, \Lambda, \chi)$ to be
$$A(\aA _{\phi}, \Lambda, \chi)=\bigoplus _{\alpha \in 
C^{\sigma}((\phi -1,\phi])}
\Lambda c_{\alpha},$$
such that the multiplication is given by 
$c_{\alpha}\ast c_{\beta}=l^{-\chi (\beta, \alpha)}c_{\alpha +\beta}.$}
\end{defi}
Note that since we assume $X$ is K3 surface or an abelian 
surface, the algebra $A(\aA _{\phi}, \Lambda, \chi)$ is a 
commutative algebra.
Let $i_{\alpha}\colon \mathfrak{Obj}^{\alpha}\aA _{\phi} \subset \mathfrak{Obj}\aA _{\phi}$ 
be the substack which parameterizes $E\in \aA _{\phi}$ of numerical type $\alpha$. 
We denote by $\Pi _{\alpha}\colon \mathfrak{Obj}^{\alpha}\aA _{\phi} 
\to \Spec \mathbb{C}$ the structure morphism. 
Given a motivic invariant $\Upsilon$ as in (\ref{mot}), we 
construct the map
 $\Theta \colon \SF (\mathfrak{Obj}\aA _{\phi}) \to 
 A(\aA _{\phi}, \Lambda, \chi)$
to be 
$$\Theta \colon f\longmapsto \sum _{\alpha \in C^{\sigma}((\phi -1,\phi])}
\Upsilon '(\Pi _{\alpha \ast}i_{\alpha}^{\ast}f) \cdot c_{\alpha}.$$
\begin{defi}\emph{\bf{\cite[Definition 3.18]{Joy4}}} \emph{
For $\alpha \in C^{\sigma}((\phi -1,\phi])$
define $\overline{\delta}^{\alpha}(\sigma)
\in A(\aA _{\phi}, \Lambda, \chi)$ to be
$$\overline{\delta}^{\alpha}(\sigma)\cneq 
\Theta (\delta ^{(\alpha, \phi(\alpha))}
(\sigma))=
I^{\alpha}(\sigma)c_{\alpha}\in A(\aA _{\phi}, \Lambda, \chi),$$
and $\overline{\epsilon}^{\alpha}(\sigma) \in A(\aA _{\phi}, \Lambda, \chi)$ 
to be 
\begin{align}\label{image}
\overline{\epsilon}^{\alpha}(\sigma)&=
\sum _{\alpha _1 +\cdots + \alpha _n =\alpha}
\frac{(-1)^{n-1}}{n}\overline{\delta}^{\alpha_1}(\sigma)\ast \cdots \ast
\overline{\delta}^{\alpha _n}(\sigma)\\
&= \frac{J^{\alpha}(\sigma)}{l-1}c_{\alpha}\in A(\aA _{\phi}, \Lambda, \chi),
\end{align}
where $\alpha _i \in C^{\sigma}(\phi(\alpha))$. }
\end{defi}
\begin{rmk}\label{motivation}\emph{
The definition of $J^{\alpha}(\sigma)$ is motivated by 
the weighted sum in the Ringel-Hall algebra. In fact 
Joyce~\cite[Theorem 8.7]{Joy3} showed that the
 following weighted sum in $\SF (\mathfrak{Obj}\aA _{\phi})$, 
$$\epsilon ^{(\alpha, \phi(\alpha))}(\sigma)=
\sum _{\alpha _1 +\cdots + \alpha _n =\alpha}
\frac{(-1)^{n-1}}{n}\delta^{(\alpha_1, \phi(\alpha _1))}(\sigma)
\ast \cdots \ast
\delta^{(\alpha _n, \phi(\alpha _n))}(\sigma), $$
with $\alpha _i \in C^{\sigma}(\phi (\alpha))$
is contained in a certain Lie subalgebra $\SF ^{\rm{ind}}(\mathfrak{Obj}
\aA _{\psi})$. Roughly it means that the stack function 
$\epsilon ^{(\alpha, \phi(\alpha))}(\sigma)$ is supported on 
indecomposable objects. Hence if $\Theta$ is a ring 
homomorphism, we have $\Theta (\epsilon ^{(\alpha, \phi(\alpha))}(\sigma))=
\overline{\epsilon}^{\alpha}(\sigma)$ and in particular one can 
define $J^{\alpha}(\sigma)|_{z=-1}\in \mathbb{Q}$
 in Example~\ref{virtual}. 
However $\Theta$ is not ring homomorphism in our case, so 
instead $J^{\alpha}(\sigma)$ is defined as 
the weighted sum in the algebra $A(\aA _{\phi}, \Lambda,
\chi)$ rather than the Ringel-Hall algebra. 
This is the motivation of the invariants explained in~\cite{Joy4}.  }
 \end{rmk}
Although the map $\Theta$ is not a ring homomorphism,  
we have the following proposition. 
\begin{prop}\label{unf}
For 
$\alpha _1, \cdots, \alpha _n \in C^{\sigma}((\phi -1,\phi])$, suppose
$\phi (\alpha _1) > \phi (\alpha _2) > \cdots > \phi (\alpha _n)$ where 
$\phi (\alpha _i)\in (\phi -1,\phi]$ the phase with respect to the 
stability function $Z$.  Then we have 
the following equality in $A(\aA _{\phi}, \Lambda, \chi)$, 
\begin{align}\label{dare}
\Theta (\delta ^{(\alpha _1, \phi (\alpha _1))}(\sigma)\ast \cdots \ast 
\delta ^{(\alpha _n, \phi (\alpha _n))}(\sigma))
=\overline{\delta}^{\alpha _1}(\sigma)
\ast \cdots \ast \overline{\delta}^{\alpha _n}(\sigma).
\end{align}
\end{prop}
\begin{proof}
This is obtained by applying~\cite[Proposition 6.20]{Joy4} 
for the abelian category $\pP ((\phi -1, \phi])$. 

\end{proof}

\subsection{Behavior of invariants in a chamber}\label{behavior}
Let us
 investigate how the invariant (\ref{yomi}) vary under a
change of stability conditions.
From here until the end of section, we fix $\alpha \in \nN (X)$.  
Let $\mathfrak{B}^{\circ} \subset \Stab ^{\ast}(X)$ be an
open subset and 
its closure $\mathfrak{B}=\overline{\mathfrak{B}}^{\circ}$ is compact. Let 
$\sS \subset D(X)$ be the set of objects, 
\begin{align}\label{mass}
\sS \cneq \{ E\in D(X) \mid 
E\mbox{ is semistable in some }\sigma '=(Z',\pP')\in \mathfrak{B} \mbox{ with }
\lvert Z'(E) \rvert \le \lvert Z'(\alpha) \rvert \}.\end{align}
Then $\sS$ has a bounded mass, thus there exists a finite 
number of codimension one submanifods $\{ \wW _{\gamma} \}_{\gamma \in \Gamma}$
which gives a wall and chamber structure on $\mathfrak{B}$.
(See Proposition~\ref{massfin}.) 
Let $\cC$ be one of the connected component, 
$$\cC \subset \mathfrak{B}\setminus \bigcup _{\gamma}\wW _{\gamma}.$$
We show the following. 
\begin{prop}\label{chamber}
Take 
$\sigma _i =(Z_i, \pP _i)\in \cC$ for $i=0,1$. 
Then we have $J^{\alpha}(\sigma _0)=J^{\alpha}(\sigma _1)$. 
\end{prop}
\begin{proof}
First assume $Z_0 (\alpha) \in \mathbb{R}_{>0}e^{i\pi \phi}$
for some $\phi \in \mathbb{R}$
 and $Z_1 (\alpha) =0$. 
Then there is no object $F\in \sS$ which is semistable 
in $\sigma _1$ and of numerical type $\alpha$. Because $\sigma _0$
and $\sigma _1$ are contained in the same chamber, there 
is no object $F\in \sS$ which is semistable in $\sigma _0$ and of numerical type 
$\alpha$. Note that if the sum (\ref{yomi}) for $\sigma _0$ is 
non-zero, there exist 
$$\alpha _1, \cdots, \alpha _n \in 
C^{\sigma _0}(\phi),$$
 such that $\prod _{i=1}^n I^{\alpha _i}(\sigma _0)$
is non-zero, and $\alpha _1 + \cdots + \alpha _n =\alpha$. 
 By the definition of $I^{\alpha _i}(\sigma _0)$, there
must be an object
$E_i \in \pP _0(\phi)$ of numerical type $\alpha _i$ for each $i$.   
Then $\oplus _{i=1}^n E_i $ is semistable in $\sigma _0$
and of numerical type $\alpha$, which is a contradiction. 
Hence in this case, one has 
$$J^{\alpha}(\sigma _0)=
J^{\alpha}(\sigma _1)=0.$$ 

Thus we may assume $Z_i (\alpha) \neq 0$ for $i=0,1$.
Again choose $\phi \in \mathbb{R}$, 
$\alpha _1, \cdots, \alpha _n \in C^{\sigma _0}(\phi)$
in the sum (\ref{yomi}) for $\sigma _0$. 
If $\prod _{i=1}^n I^{\alpha _i}(\sigma _0) \neq 0$, 
then there exist $E _i \in \pP _0(\phi)$ of numerical type $\alpha _i$. 
We have 
$$\lvert Z_0 (E_i) \rvert =\lvert Z_0 (\alpha _i) \rvert 
\le \lvert Z_{0}(\alpha) \rvert .$$
Thus we have $E_i \in \sS$. Note that for each $i$ and $j$, 
the values $Z_{0}(\alpha _i)$ and $Z_{0}(\alpha _j)$ are 
proportional. Furthermore 
we have $\sigma _0 \notin \wW _{\gamma}$ for 
any $\gamma$.  By the construction of $\wW _{\gamma}$ in 
Proposition~\ref{massfin}, 
this implies $\alpha _i$ and $\alpha _j$ 
must be proportional in $\nN (X)$, hence $\alpha _i$ is proportional 
to $\alpha$. 
Choose a path 
$$\lambda \colon [0,1] \lr \cC,$$
such that $\lambda (0)=\sigma _0$ and $\lambda (1)=\sigma _1$.
We denote $\lambda (t)=\sigma _t =(Z_t, \pP _t)$. 
For an arbitrary $E_i 
\in \pP _0(\phi)$ of numerical type $\alpha _i$, we have 
$E_i \in \sS$, thus $E_i$ is also semistable 
in $\sigma _{t}$. Hence $Z_t (\alpha _i)\neq 0$ for $t\in [0,1]$, and
the phase of $E_i$ in $\sigma _{t}$ is uniquely determined 
independent of a choice of $E_i \in \pP _0 (\phi)$. 
Thus there is $\phi ' \in \mathbb{R}$ which satisfies 
$Z_{1}(\alpha _i)\in \mathbb{R}_{>0}e^{i\pi \phi '}$ such that 
$$\mM ^{(\alpha _i, \phi)}(\sigma _0) \subset 
\mM ^{(\alpha _i, \phi ')}
(\sigma _{1}).$$
By the converse argument, we obtain $\mM ^{(\alpha _i, \phi)}
(\sigma _0) =
\mM ^{(\alpha _i, \phi ')}
(\sigma _{1})$, thus 
$I^{\alpha _i}(\sigma _0)=I^{\alpha _i}(\sigma _{1})$. 
Because $\alpha _i$ is proportional to $\alpha$, we have 
$Z_{1}(\alpha)\in \mathbb{R}_{>0}e^{i\pi \phi '}$ and 
$\alpha _i \in C^{\sigma _1}(\phi ')$. 
Hence the sum (\ref{yomi}) for $J^{\alpha}(\sigma _0)$ and 
$J^{\alpha}(\sigma _{1})$ are identified. 
\end{proof}

\subsection{Behavior of invariants near a wall}

Next we investigate the behavior of the
invariants near a wall $\wW _{\lambda}$. 
Here we use the same notation as in (\ref{behavior}). 
Take $0<\varepsilon <1/6$ and 
$\sigma _i=(Z_i, \pP _i)\in \Stab ^{\ast}(X)$
 for $i=0,1$. We assume the following. 
\begin{itemize}
\item    
$\sigma _0$ is algebraic, contained in $\cC$, 
 and
$\pP ((\psi -1, \psi])$ satisfies the generic flatness
for any $\psi \in \iI$.  (See Remark~\ref{remark}.)
\item $\sigma _1 \in \wW _{\gamma}\cap \overline{\cC}$
for some $\gamma$
and $\sigma _0 \in B_{\varepsilon}(\sigma _1)$. (See (\ref{opbasis}).)
\end{itemize}
First we give the following lemma. 
\begin{lem}\label{zero}
Assume $Z_1 (\alpha)=0$. Then we have 
$$J^{\alpha}(\sigma _0)=J^{\alpha}(\sigma _1)=0.$$
\end{lem}
\begin{proof}
By the definition, we have $J^{\alpha}(\sigma _1)=0$. 
Assume that there is an object $E\in D(X)$, semistable in $\sigma _0$ of 
numerical type $\alpha$. Then $E\in \sS$, hence $E$ is 
semistable for arbitrary 
$\sigma _0 '\in \cC$. By the comment in~\cite{Brs1} 
after~\cite[Proposition 8.1]{Brs1}, $E$ is also semistable in $\sigma _0$
hence a contradiction. Therefore there is no semistable object in $\sigma _0$
of numerical type $\alpha$, and this implies $J^{\alpha}(\sigma _0)=0$. 
\end{proof}
By Lemma~\ref{zero}, it is enough to consider the case of 
$Z_1 (\alpha)\neq 0$.
Choose $\phi \in \mathbb{R}$ which satisfies 
$Z_1 (\alpha)\in \mathbb{R}_{>0}e^{i\pi \phi}$. 
Note that for any 
$\beta \in C^{\sigma _0}((\phi -\varepsilon, \phi + \varepsilon))$,
we can define the phases 
$$\phi _0(\beta)\in (\phi -\varepsilon, \phi +\varepsilon), \quad 
\phi _{1}(\beta) \in (\phi -2\varepsilon, \phi +2\varepsilon),$$
with respect to the stability functions $Z_0$, $Z_{1}$ respectively. 
We fix $\psi \in \iI$ which satisfies 
$$\pP _1 (\phi) \subset \pP _0 ((\phi -\varepsilon, \phi +\varepsilon))
\subset \aA _{\psi}=\pP _0 ((\psi-1, \psi]).$$
We
consider the $\mathbb{Q}$-algebra $(\SF (\mathfrak{Obj}\aA _{\psi}), \ast)$. 
\begin{lem}\label{sosite}
We have the following equality in $\SF (\mathfrak{Obj}\aA _{\psi})$, 
\begin{align}\label{desita}\delta ^{(\alpha, \phi)}(\sigma _1)=
\sum _{\alpha _1 + \cdots + \alpha _n = \alpha}
\delta ^{(\alpha _1, \phi _0(\alpha _1))}(\sigma _0)
\ast \cdots \ast \delta ^{(\alpha_n, \phi_0 (\alpha _n))}
(\sigma _0),\end{align}
where $\alpha _i \in C^{\sigma _0}((\phi -\varepsilon, \phi +\varepsilon))$, 
and $\{ \alpha _i \}_{1\le i \le n}$
satisfy
\begin{align}\label{modoru}
\phi _0 (\alpha _1) > \cdots > \phi _0 (\alpha _n), \quad 
\phi _1 (\alpha _1)= \cdots = \phi _1 (\alpha _n)=\phi.\end{align}
\end{lem}
\begin{proof}
We show the following decomposition, 
\begin{align}\label{decompos}
\mM ^{(\alpha, \phi)}(\sigma _1)=
\coprod _{\alpha _1 +\cdots + \alpha _n =\alpha}
\Pi _{n\ast} \mM (\{ \alpha _i\}_{1\le i\le n}, \aA _{\psi}, \sigma _0),
\end{align}
where $\alpha _i \in C^{\sigma _0}((\phi -\varepsilon, \phi +\varepsilon))$
and they satisfy (\ref{modoru}). 
First note that any object $E\in \pP _1 (\phi)$ 
of numerical type $\alpha$, a $\mathbb{C}$-valued point of the LHS 
of (\ref{decompos}), 
has the unique filtration
$$0=E_0 \subset E_1 \subset \cdots \subset E_n =E,$$
in $\aA _{\psi}$ such that $F_i =E_i /E_{i-1}$
is semistable in $\sigma _0$,   
of numerical type $\alpha _i \in C^{\sigma _0}
((\phi -\varepsilon, \phi +\varepsilon))$, and they
satisfy $\phi _0 (\alpha _1)> \cdots > \phi _0 (\alpha _n)$.  
Since $0<\varepsilon <1/6$, we have 
$$\lvert Z_0(F_i) \rvert \le \lvert Z_0 (E) \rvert =
\lvert Z_0 (\alpha)\rvert,$$
thus $F_i \in \sS$. 
Therefore $F_i$ is semistable for any $\sigma _0' \in \cC$, hence 
it is also semistable in $\sigma _1$. 
The condition $\sigma _1 \in \overline{\cC}$
implies that $\phi _1(\alpha _1) \ge \cdots \ge \phi _1 (\alpha _n)$. 
Because $E$ is semistable in $\sigma _1$, 
we must have $\phi _1 (\alpha _1)= \cdots = \phi _1 (\alpha _n)=\phi$.
This means $E$ is a $\mathbb{C}$-valued point of
the RHS of (\ref{decompos}). 
 
Conversely take $\alpha _1, \cdots, \alpha _n
 \in C^{\sigma _0}((\phi -\varepsilon, \phi +\varepsilon))$ which 
 satisfy $\alpha _1 + \cdots + \alpha _n =\alpha$ and (\ref{modoru}). 
 Suppose for $E\in \aA _{\psi}$, there is a filtration 
$$0=E_0 \subset E_1 \subset \cdots \subset E_n =E,$$ in 
$\aA _{\psi}$ such that $F_i =E_i/E_{i-1}$ is semistable 
in $\sigma _0$ and of numerical type $\alpha _i$, i.e. $E$ is a
$\mathbb{C}$-valued point of the RHS of (\ref{decompos}). 
 Again we have 
$\lvert Z_0 (F_i) \rvert \le \lvert Z_0 (\alpha) \rvert$ thus 
$F_i \in \sS$. Hence $F_i$ is also semistable in $\sigma _1$, 
and (\ref{modoru}) implies $E$ is also semistable in $\sigma _1$. 
This implies $E$ is a $\mathbb{C}$-valued point 
of the LHS of (\ref{decompos}). 

Now we have shown (\ref{decompos}) at the level 
of $\mathbb{C}$-valued points. Finally 
we have the isomorphism of the stabilizers, 
$$\Aut (E_1\subset \cdots \subset E_n) \cong \Aut (E_n).$$
Hence we have the decomposition (\ref{decompos}). 
Using (\ref{kowareta}), 
the formula (\ref{desita}) follows. 
(Also see the proof of~\cite[Theorem 5.11]{Joy4}.)
\end{proof}
Next we compare $\overline{\epsilon} ^{\alpha}(\sigma _i)
 \in A(\aA _{\psi}, \Lambda, \chi)$ near a wall. 
 Following~\cite[Definition 4.2]{Joy4}, we introduce the following 
 combinatorial values.  
\begin{defi}\label{comb}
\emph{ For $\alpha _1, \cdots, \alpha _n \in C^{\sigma _0}
((\phi -\varepsilon, \phi +\varepsilon))$, 
consider the following two conditions.}

\emph{(a) $\phi _0(\alpha _i) \le \phi _0 (\alpha _{i+1})$, and 
$\phi _{1}(\alpha _1 +\cdots + \alpha _i)> \phi _{1}(\alpha _{i+1}
+ \cdots \alpha _n)$. }

\emph{(b) $\phi _0(\alpha _i) > \phi _0 (\alpha _{i+1})$, and 
$\phi _{1}(\alpha _1 +\cdots + \alpha _i)\le \phi _{1}(\alpha _{i+1}
+ \cdots \alpha _n)$. }

\emph{If for all $i=1, \cdots, n-1$ one 
of the above two conditions is satisfied, 
then define
$$S(\{ \alpha _i \}_{1\le i\le n}, \sigma _0, \sigma _1)=(-1)^{r},$$
where $r$ is the number of $i=1,\cdots,n-1$ satisfying $(a)$.  
Otherwise define $S(\{ \alpha _i \}_{1\le i\le n}, \sigma _0, \sigma _1)=0$. }
\end{defi} 
The values $S(\{ \alpha _i \}_{1\le i\le n}, \sigma _0, \sigma _1)$
give the transformation coefficients of the invariants. 
\begin{lem}\label{coe}
We have the following equality in $A(\aA _{\phi}, \chi, \Lambda)$, 
\begin{align}\label{eq}
\overline{\delta}^{\alpha}(\sigma _1)=
\sum _{\alpha _1 + \cdots +\alpha _n 
=\alpha}S(\{ \alpha _i \}_{1\le i\le n}, \sigma _0, 
\sigma _1) \overline{\delta}^{\alpha _1}(\sigma _0) \ast \cdots \ast 
\overline{\delta}^{\alpha _n}(\sigma _0),\end{align}
where $\alpha _i 
\in C^{\sigma _0}((\phi -\varepsilon, \phi +\varepsilon ))$.
\end{lem}
\begin{proof}
Applying $\Theta$ to (\ref{desita}) and using Proposition~\ref{unf}, 
we have 
\begin{align}\label{kokomade}
\overline{\delta}^{\alpha}(\sigma _1)=
\sum _{\alpha _1 + \cdots +\alpha _n =\alpha}\overline{\delta}^{\alpha _1}
(\sigma _0) \ast \cdots \ast 
\overline{\delta}^{\alpha _n}(\sigma _0),
\end{align} 
where $\alpha _1, \cdots, \alpha _n \in C^{\sigma _0}((\phi -\varepsilon, 
\phi +\varepsilon))$ satisfy (\ref{modoru}). 
Therefore it is enough to check that the right hand sides of 
(\ref{kokomade}) and (\ref{eq}) are equal. 
Suppose that 
$\alpha _1, \cdots, \alpha _n$ in (\ref{eq}) satisfy
\begin{align}\label{mune}
S(\{ \alpha _i \}_{1\le i\le n}, \sigma _0, \sigma _1) \overline{\delta}^{\alpha _1}(\sigma _0) \ast \cdots \ast \overline{\delta}^{\alpha _n}(\sigma _0) \neq
\emptyset.\end{align}
Then for each $i$, there exists $E_i \in \pP _0 (\phi _0(\alpha _i))$ 
which is of numerical type $\alpha _i$. 
Since $0< \varepsilon <1/6$, we have 
$\lvert Z_0 (E_i) \rvert \le \lvert Z_0 (\alpha) \rvert$.
Thus $E_i \in \sS$, and by the construction of walls $\wW _{\lambda}$
in Proposition~\ref{massfin}, 
we have the following: 
\begin{align}\label{kite}
\mbox{For } i,j, \mbox{ if } \phi _0 (\alpha _i) \ge 
\phi _0 (\alpha _j) \mbox{ then }
\phi _1 (\alpha _i) \ge \phi _1(\alpha _j).\end{align} 
Thus the coefficient $S(\{ \alpha _i \}_{1\le i\le n}, \sigma _0, \sigma _1)$
for which $\alpha _1, \cdots, \alpha _n$ satisfy (\ref{mune}) 
is calculated by the property (\ref{kite})
 in a purely combinatorial way. 
 It is computed in~\cite[5.2]{Joy4} and the result is 
$$S(\{ \alpha _i \}_{1\le i\le n}, \sigma _0, \sigma _1)=\left\{
\begin{array}{ll}1 & \mbox{ if }\alpha _1, \cdots, \alpha _n 
\mbox{ satisfy (\ref{modoru})}, \\
0 & \mbox{ otherwise }.
\end{array}\right.$$
\end{proof}

\begin{rmk}\label{rara}\emph{
Take $\alpha '\in C^{\sigma _1}(\phi)$ with
 $\lvert Z_1(\alpha ') \rvert \le 
\lvert Z_1 (\alpha) \rvert$. Then the above proof shows that 
the same formula (\ref{eq}) replaced $\alpha$ by $\alpha '$ also holds. }
\end{rmk}
Now we can compare $J^{\alpha}(\sigma _0)$ and $J^{\alpha}(\sigma _1)$. 
\begin{prop}\label{kiite}
We have $\overline{\epsilon}^{\alpha}(\sigma _0)=\overline{\epsilon}^{\alpha}
(\sigma _1)$ in $A(\aA _{\psi}, \chi, \Lambda)$.
 Thus we have $J^{\alpha}(\sigma _0)=J^{\alpha}(\sigma _1)$ in 
 $\Lambda$. 
\end{prop}
The proof relies on the 
the combinatory result in~\cite{Joy4}.
So before beginning the proof, we show the simplest case
for the reader's convenience.  
Suppose the following decomposition is unique. 
$$\alpha =\alpha _1 + \alpha _2, \quad \alpha _j \in 
C^{\sigma _0}((\phi -\varepsilon, \phi +\varepsilon)).$$
 Furthermore for such  $\alpha _j$, we assume
  $\phi _1 (\alpha _1)=\phi _1 (\alpha _2)$,  
  $\phi _0 (\alpha _1)> \phi _0 (\alpha _2)$. 
In this case we have 
\begin{align*}
\overline{\epsilon}^{\alpha}(\sigma _1)&= \overline{\delta}^{\alpha}(\sigma _1)
-\frac{1}{2}\overline{\delta}^{\alpha _1}(\sigma _1)\ast 
\overline{\delta}^{\alpha _2}(\sigma _1)
-\frac{1}{2}\overline{\delta}^{\alpha _2}(\sigma _1)\ast 
\overline{\delta}^{\alpha _1}(\sigma _1), \\
\overline{\delta}^{\alpha}(\sigma _1) &= \overline{\delta}^{\alpha}(\sigma _0)
+ \overline{\delta}^{\alpha _1}(\sigma _0)\ast 
\overline{\delta}^{\alpha _2}(\sigma _0).
\end{align*}
Also we have $\overline{\epsilon}^{\alpha}(\sigma _0)
=\overline{\delta}^{\alpha}(\sigma _0)$, and 
$\overline{\delta}^{\alpha _j}(\sigma _i)
=\overline{\epsilon}^{\alpha _j}(\sigma _i)$. Thus we have 
\begin{align*}
\overline{\epsilon}^{\alpha}(\sigma _1)&= 
\overline{\epsilon}^{\alpha}(\sigma _0) + \frac{1}{2}
[\overline{\epsilon}^{\alpha _1}(\sigma _0),
 \overline{\epsilon}^{\alpha _2}(\sigma _0)], \\
&= \overline{\epsilon}^{\alpha}(\sigma _0).\end{align*}
Here we used the fact that $A(\aA _{\psi}, \chi, \Lambda)$ is commutative 
in our case. Now we give the proof of Proposition~\ref{kiite}. 

\begin{proof}
Using the proof of~\cite[Theorem 7.7]{Joy3}
in the algebra $A(\aA _{\psi}, \Lambda, \chi)$
(also see~\cite[Definition 6.22]{Joy4}), 
 we may write $\overline{\delta} ^{\alpha}(\sigma _0)$
as follows,  
\begin{align}\label{may}
\overline{\delta}^{\alpha}(\sigma _0)
=\sum _{\alpha _1 +\cdots \alpha _n =\alpha}
\frac{1}{n!}\overline{\epsilon}^{\alpha _1}(\sigma _0)\ast \cdots 
\ast \overline{\epsilon}^{\alpha _n}(\sigma _0),\end{align}
where $\alpha _i \in C^{\sigma _0}(\phi _0 (\alpha))$. 
We can rewrite (\ref{may}) as the same formula of (\ref{may}) 
and $\alpha _i \in C^{\sigma _0}((\phi -\varepsilon, \phi +\varepsilon))$
with $\phi _0 (\alpha _i)=\phi _0(\alpha)$. 
Then substituting (\ref{eq}), (\ref{may}) to the definition of 
$\overline{\epsilon}^{\alpha}(\sigma _1)$ in (\ref{image}), 
(also noting Remark~\ref{rara})
we can write $\overline{\epsilon}^{\alpha}(\sigma _1)$
in the following formula,
\begin{align}\label{kaze}
\overline{\epsilon}^{\alpha}(\sigma _1)=\sum _{\alpha _1 +\cdots + 
\alpha _n =\alpha}
U(\{ \alpha _i\}_{1\le i\le n}, \sigma _0, \sigma _1)\overline{\epsilon}^{\alpha _1}
(\sigma _0)
\ast \cdots \ast \overline{\epsilon}^{\alpha _n}(\sigma _0),\end{align}
where $\alpha _i \in C^{\sigma _0}((\phi -\varepsilon, \phi +\varepsilon))$.
 Here one can consult the explicit description of 
 $$U(\{\alpha _i\}_{1\le i\le n}, \sigma _0, \sigma _1 )\in \mathbb{Q},$$
 in~\cite[Definition 4.4]{Joy4}, after replacing $C(\aA)$ 
 in \textit{loc.cite.}
 by $C^{\sigma _0}((\phi -\varepsilon, \phi +\varepsilon))$.
 In~\cite[Theorem 5.2]{Joy4}, it is proved that
  using~\cite[Theorem 5.4]{Joy4} the formula (\ref{kaze}) is 
  written as 
 $$\overline{\epsilon}^{\alpha}(\sigma _1)=
 \overline{\epsilon}^{\alpha}(\sigma _0) + [\mbox{ multiple commutators of }
 \overline{\epsilon}^{\alpha _i}(\sigma _0)].$$
In our case $A(\aA _{\psi}, \Lambda, \chi)$ is commutative since we assume  
 $X$ is a K3 surface or an abelian surface. Hence we have 
 $\overline{\epsilon}^{\alpha}(\sigma _1)=
 \overline{\epsilon}^{\alpha}(\sigma _0)$.
 \end{proof}

Now we can show the following. 
\begin{thm}\label{inde}
For $\sigma =(Z, \pP)\in \Stab ^{\ast}(X)$ and $\alpha \in\nN (X)$, 
the invariant $J^{\alpha}(\sigma)\in \Lambda$ does not 
depend on a choice of $\sigma$. 
\end{thm}

\begin{proof}
Take $\sigma =(Z, \pP)\in \Stab ^{\ast}(X)$ and 
$\tau =(W, \qQ)\in \Stab ^{\ast}(X)$.
Let
$\lambda$ be a path 
$$\lambda \colon [0,1] \lr  \Stab ^{\ast}(X),$$
such that $\lambda (0)=\sigma$ and $\lambda (1)=\tau$. 
We take a connected open set $\mathfrak{B}^{\circ}\subset \Stab ^{\ast}(X)$
which contains $\lambda ([0,1])$ and its closure
$\mathfrak{B}=\overline{\mathfrak{B}}^{\circ}$ is compact. 
We consider the set of objects $\sS$ as in (\ref{mass}) 
and the associated walls $\{ \wW _{\gamma} \}_{\gamma \in \Gamma}$. 
We denote $\lambda (t)=\sigma _t =(Z_t, \pP _t)$. 
We may assume that the set of points $K\subset [0,1]$ on which 
$\sigma _t$ is algebraic and $\pP _t ((\psi -1, \psi])$ satisfies 
the generic flatness for any $\psi \in \iI$
is dense in $[0,1]$. 
Take $s_0, s_1, s_2, \cdots, s_N, s_{N+1} \in [0,1]$ 
and $t_i ^{\pm}\in (s_i, s_{i+1})\cap K$ such that 
\begin{itemize}
\item For $1\le i\le N$, $s_i \in \wW _{\gamma}$ for some $\wW _{\gamma}$,
and $s_0=0$, $s_{N+1}=1$. 
\item For any $t\in (s_i, s_{i+1})$, 
we have $\lambda (t) \notin \wW _{\gamma}$ for any $\gamma$. 
\item $\sigma _{t_i ^{+}} \in B_{\varepsilon}(\sigma _{s_{i+1}})$, 
$\sigma _{t_{i}^{-}}\in 
B_{\varepsilon}(\sigma _{s_{i}})$ with $0<\varepsilon <1/6$.  
\end{itemize}
By Proposition~\ref{chamber} and Proposition~\ref{kiite}, we have 
$$J^{\alpha}(\sigma _{s_i})=J^{\alpha}(\sigma _{t_i ^{-}})
= J^{\alpha}(\sigma _{t_i ^{+}}) =J^{\alpha}(\sigma _{s_{i+1}}),$$
for each $i$. Thus $J^{\alpha}(\sigma)=J^{\alpha}(\tau)$ follows. 
\end{proof}

By Theorem~\ref{inde}, the following definition is well defined. 
\begin{defi}\emph{
For $\alpha \in\nN (X)$, we define $J^{\alpha}\in \Lambda$ 
to be $J^{\alpha}(\sigma)$ for some $\sigma \in \Stab ^{\ast}(X)$. }
\end{defi}
Let $\Auteq ^{\ast}D(X)$ be the subgroup of $\Phi \in \Auteq D(X)$ 
which preserves the connected component $\Stab ^{\ast}(X)$. Also 
for $\Phi \in \Auteq ^{\ast}D(X)$, we denote 
$$\Phi _{\ast}\colon \nN (X) \lr \nN (X),$$
the induced automorphism. 
We have the following corollary of Theorem~\ref{inde}. 
\begin{cor}\label{yume}
For $\Phi \in \Auteq ^{\ast}D(X)$, one has 
$J^{\alpha}=J^{\Phi _{\ast}\alpha}$. 
\end{cor}
\begin{proof}
We have 
$$J^{\alpha}=J^{\alpha}(\sigma)= J^{\Phi _{\ast}\alpha}(\Phi (\sigma))
= J^{\Phi _{\ast}\alpha}.$$
\end{proof}

\section{Comparison of invariants which count semistable objects and
semistable sheaves}
\label{compare}
In this section we compare $J^{\alpha}$ and $\hat{J}^{\alpha}$, 
where $\hat{J}^{\alpha}$ is a
counting invariant of semistable sheaves introduced in~\cite{Joy4}. 
\subsection{Counting invariants of semistable sheaves}
Let 
$\Lambda$ be a $\mathbb{Q}$-algebra and
$\Upsilon \colon K(\Var) \to \Lambda$
be a motivic invariant as in (\ref{mot}). 
We denote by $C(X)$ the image, 
$$C(X) = \im (\Coh (X) \lr \nN (X)).$$
For $\alpha \in C(X)$, we recall 
the definition of $\hat{J}^{\alpha} \in \Lambda$ 
introduced in~\cite{Joy4}. 
Let $\omega$ be an ample divisor on $X$. We consider the 
moduli stack, 
$$\hat{\mM}^{\alpha}(\omega)\subset \mM,$$
which is the stack of $\omega$-Gieseker semistable sheaves of 
numerical type $\alpha$.  Let 
$\hat{\delta}^{\alpha}(\omega) \in \SF (\mM)$ be the associated 
stack function. 
We consider the map $\Upsilon '\circ \Pi _{\ast}\colon 
\SF (\mM) \to \Lambda$ as in (\ref{signal}). 
\begin{defi}\emph{\bf{\cite[Definition 6.1, Definition 6.22]{Joy4}}}
\emph{ We define $\hat{I}^{\alpha}(\omega) \in \Lambda$ to be
 $$\hat{I}^{\alpha}(\omega)= \Upsilon ' \circ \Pi _{\ast}\hat{\delta}^{\alpha}
 (\omega) \in \Lambda,$$
 and $\hat{J}^{\alpha}(\omega) \in \Lambda$ to be 
 \begin{align}\label{countsheaf}
 \hat{J}^{\alpha}(\omega)=
\sum _{\alpha _1 + \cdots +\alpha _n=\alpha}
l^{-\sum _{j>i}\chi (\alpha _j, \alpha _i)}
\frac{(-1)^{n-1}(l-1)}{n}\prod _{i=1}^n \hat{I}^{\alpha _i}
(\omega) \in \Lambda,
\end{align}
where $\alpha _i \in C(X)$ satisfies $P(\alpha _i, \omega, n)= P(\alpha, 
\omega, n)$. }
\end{defi}
Joyce~\cite{Joy4} showed the following. 
\begin{thm}\emph{\bf{\cite[Theorem 6.24]{Joy4}}}\label{indamp} 
The invariant $\hat{J}^{\alpha}(\omega)\in \Lambda$
does not depend on a choice of an ample divisor $\omega$. 
\end{thm}
For $\alpha \in C(X)$, we define $\hat{J}^{\alpha}\in \Lambda$
to be $\hat{J}^{\alpha}(\omega)$ for some ample divisor $\omega$, 
which is well defined by Theorem~\ref{indamp}.

\subsection{Comparison of $J^{\alpha}$ and $\hat{J}^{\alpha}$}
Here we compare $J^{\alpha}$ and $\hat{J}^{\alpha}$ for $\alpha \in C(X)$. 
Let us take an ample divisor $\omega$ and $k\in \mathbb{Q}_{\ge 1}$.  
We use the following notation, 
 $$Z_{k\omega} =Z_{(0,k\omega)}, \quad 
 \aA _{\omega}=\aA _{(0,\omega)}=\aA _{(0, k\omega)}, 
 \quad \sigma _k =(Z_{k\omega}, \aA _{\omega}).$$
The idea is to compare the following two values, 
$$J^{\alpha}(\sigma _{k})\in \Lambda, \quad
\hat{J}^{\alpha}(\omega)\in \Lambda,$$
in the limit $k\to \infty$. 
In~\cite[Proposition 14.2]{Brs2}, Bridgeland
proved that (putting a certain assumption on a numerical class), 
an object $E\in D(X)$ is semistable in $\sigma _k$
for all $k\gg 0$ if and only if $E$ is $\omega$-Gieseker semistable. 
This is what string theory predicts that BPS branes in the limit 
$k\to \infty$ are in fact Gieseker stable sheaves. What we actually have 
to prove is that we can choose $k>0$ uniformly so that it works 
for any semistable objects. First we give the following lemma. 
\begin{lem}\label{line}
For an ample line bundle $\lL \in \Pic (X)$, one has 
\begin{align}\label{gensou}
J^{\alpha}=J^{\alpha \otimes \lL}, \quad 
\hat{J}^{\alpha}=\hat{J}^{\alpha \otimes \lL}.
\end{align}

\end{lem}
\begin{proof}
Note that tensoring $\lL$ gives an 
autoequivalence 
$\otimes \lL \in \Auteq ^{\ast}D(X)$. 
Thus $J^{\alpha}=J^{\alpha \otimes \lL}$ follows from Corollary~\ref{yume}. 
Next let $\omega =c_1 (\lL)$. Then by Theorem~\ref{indamp}
we have $\hat{J}^{\alpha}=\hat{J}^{\alpha}(\omega)$. 
The equality $\hat{J}^{\alpha}(\omega)=\hat{J}^{\alpha \otimes \lL}(\omega)$ 
follows easily from the fact that for any $\omega$-Gieseker semistable sheaf
$E$ of numerical type $\alpha$, $E\otimes \lL$ is also $\omega$-Gieseker semistable and 
it is of numerical type $\alpha \otimes \lL$. 
\end{proof}
For $\alpha \in C(X)$ we denote $v(\alpha)=(r,l,s)$. 
We show the following proposition. 
\begin{prop}\label{semi}
Suppose $\omega \cdot l >0$ or $r=l=0$, and 
choose $0<\phi _k \le 1$ which satisfies 
$Z_{k\omega} (\alpha)\in \mathbb{R}_{>0}e^{i\pi \phi _k}$. 
Then
there exists $N>0$ such that for all $k\ge N$ and 
$\alpha '$ which satisfies 
\begin{align}\label{ei}
\alpha ' \in C^{\sigma _k}(\phi _k) \mbox{ with }
\lvert \Imm Z_{\omega}(\alpha ') \rvert \le 
\lvert \Imm Z_{\omega}(\alpha) \rvert, 
\end{align}
any $E\in M^{(\alpha ', \phi _k)}(\sigma _k)$ is a $\omega$-Gieseker 
semistable sheaf. 

\end{prop}
\begin{proof}\

\begin{ssstep}\end{ssstep}
First by Lemma~\ref{useful} (i), (ii), the set of $\alpha ' \in \nN (X)$
which satisfies (\ref{ei}) is a finite set for a fixed $\alpha$. 
When $r=l=0$, any object $E\in \aA _{\omega}$ of
numerical type $\alpha$ is a zero dimensional sheaf, so the result is 
obvious. Thus we may assume $\omega \cdot l >0$. 
In this case $\phi _k$ goes to zero for $k\to \infty$
when $r>0$ and goes to $1/2$ when $r=0$. 
Thus there is $N>0$ so that $\phi _k \le 3/4$ for all $k\ge N$.  
Take $E\in M^{(\alpha ', \phi _k)}(\sigma _k)$, and $\alpha '$
satisfies (\ref{ei}). Then we have 
$$\phi _k (H^{-1}(E)[1]) \le \phi _k \le \frac{3}{4}.$$
Thus the map on $\cup _{k\ge N, \alpha '}M^{(\alpha ', \phi _k)}(\sigma _k)$,
\begin{align*}E 
& \longmapsto \frac{\Ree Z_{k\omega}(H^{-1}(E)[1])}
{\Imm Z_{k\omega}(H^{-1}(E)[1])} \\
& = \frac{1}{k}\cdot \frac{\Ree Z_{k\omega}(H^{-1}(E)[1])}{\Imm Z_{\omega}
(H^{-1}(E)[1])},\end{align*}
is bounded below. 
Note that $E$ is contained in $M^{\alpha}(0, \omega)$, and the map
$E\mapsto \Imm Z_{\omega}(E)$ 
on $M^{\alpha}(0, \omega)$ is bounded by Lemma~\ref{con1}. 
Therefore the map on 
$\cup _{k\ge N, \alpha '}M^{(\alpha ', \phi _k)}(\sigma _k)$, 
$$E \longmapsto \frac{1}{k} \Ree Z_{k\omega}(H^{-1}(E)[1]),$$
is bounded below. 
Thus using Lemma~\ref{fina}, we may assume that 
$$E\longmapsto \Ree Z_{\omega}(H^{-1}(E)[1]),$$
is bounded below
on $\cup _{k\ge N, \alpha '}M^{(\alpha ', \phi _k)}(\sigma _k)$.
 Then one can apply Lemma~\ref{asu}
and the set 
$$\{ v(H^{-1}(E)[1]) \in \NS ^{\ast}(X) \mid E\in \bigcup _{k\ge N, \alpha '}
M^{(\alpha ', \phi _k)}(\sigma _k)\},$$
is a finite set. Let us denote the above set
$v_1, \cdots, v_n$. Then $\lim _{k\to \infty}\phi _k(v_i)=1$, 
so by replacing $N$ if necessary, we have $\phi _k(v_i)>3/4$ 
for any $k\ge N$ and $1\le i\le n$. This
implies that for $k\ge N$ any object $E\in M^{(\alpha ', \phi _k)}(\sigma _k)$
satisfies $H^{-1}(E)=0$, so $E$ is a sheaf.

\begin{ssstep}\end{ssstep}
Next we show that any 
$E\in \cup _{k\ge N, \alpha '}M^{(\alpha ', \phi _k)}(\sigma _k)$
is a $\omega$-Gieseker semistable sheaf, by replacing $N$ if necessary. 
Assume that $E$ is not $\omega$-Gieseker semistable, and 
let $T$ be the $\omega$-Gieseker semistable factor of $E$ 
of the smallest reduced Hilbert polynomial. We denote 
$$v(E)=(r', l', s'), \quad
v(T)=(r'',l'',s'').$$
Note that $r'=0$ is equivalent to $r''=0$, and in this 
case $E$ must be $\omega$-semistable by Lemma~\ref{useful} (iii).  
Thus we may assume $r'>0$, $r''>0$. Note that 
in this case $\phi _k$ goes to 0 for $k \to \infty$. 
Since the map $E\to T$ is a surjection in $\aA _\omega$, 
and $E$ is $\sigma _k$-semistable for some $k\ge N$, one has 
$\phi _k(E) \le \phi _k(T)$. Thus we have 
\begin{align}\label{setuna}
\frac{\Ree Z_{k\omega}(E)}{\Imm Z_{k\omega}(E)} \ge
\frac{\Ree Z_{k\omega}(T)}{\Imm Z_{k\omega}(T)}.\end{align}
Explicitly (\ref{setuna}) is equivalent to 
\begin{align}\label{setuna2}
\frac{\omega\cdot l''}{\omega \cdot l'}\left(-s' +\frac{1}{2}k^2 \omega ^2 r'
\right) \ge -s''+ \frac{1}{2}k^2 \omega ^2 r''.
\end{align}
Also we note that 
\begin{align}\label{setuna3}
0< r'' \le r', \quad 0< \omega \cdot l'' \le \omega \cdot l', 
\quad l^{''2}-2r''s''\ge -2.
\end{align}
Here the third equality comes from Lemma~\ref{useful} (i). Then 
(\ref{setuna2}) and (\ref{setuna3}) imply that 
the set 
\begin{align}\label{kirei}
\{ v(T) \in \NS ^{\ast}(X) \mid E\in \cup _{k\ge N, \alpha '}
M^{(\alpha ', \phi _k)}(\sigma _k)\},\end{align}
is a finite set. 
Applying the same argument for other torsion free 
Gieseker-semistable factors, we deduce that the set 
$$\{ v(E_{\rm{fr}}) \in \NS ^{\ast}(X) \mid E\in \cup _{k\ge N, \alpha '}
M^{(\alpha ', \phi _k)}(\sigma _k)\},$$
is a finite set. It follows that the set
$$\{ v(E_{\rm{tor}}) \in \NS ^{\ast}(X) \mid E \in \cup _{k\ge N, \alpha '}
M^{(\alpha ', \phi _k)}(\sigma _k)\},$$
is also a finite set, say $v_1', \cdots, v_m'$. 
Since $\phi _k(v_i')$ goes to $1/2$ for each $i$, 
we have $\phi _k(v_i ') > \phi _k$ for all $1\le i\le m$ and $k\ge N$, 
after replacing $N$ if necessary. Thus for such $N$ and $k\ge N$,  
if we take $E\in M^{(\alpha ', \phi _k)}(\sigma _k)$, then
 $E$ must be a torsion free sheaf. 
By the definition of $T$, one has 
$$\mu _{\omega}(E)>\mu _{\omega}(T) \quad \mbox{ or } \quad
\mu _{\omega}(E)=\mu _{\omega}(T), \frac{s'}{r'}>\frac{s''}{r''}.$$
We have 
$$\frac{Z_{k\omega}(E)}{r'}-\frac{Z_{k\omega}(T)}{r''}
=-\left( 
\frac{s'}{r'}-\frac{s''}{r''}\right)+ik (\mu _{\omega}(E)-\mu _{\omega}(T)).$$
So after replacing $N$ we have $\phi _k(E)>\phi _k (T)$ 
for $k\ge N$. Such $N$ is determined by only a
numerical class of $T$. Thus the finiteness of (\ref{kirei})
implies that one can take $N$ uniformly so that 
$\phi _k(E)>\phi _k(T)$ for all $E\in M^{(\alpha ', \phi _k)}(\sigma _k)$
and $k\ge N$. This contradicts that $E$ is $\sigma _k$-semistable, so
$E$ must be $\omega$-Gieseker semistable. 
\end{proof}
Next we check the following. 

\begin{lem}\label{koware}
Suppose $\omega \cdot l>0$ or $r=l=0$. Then
there is $N>0$ so that for $k\ge N$
and $\alpha '\in C(X)$ which satisfies
\begin{align}\label{hirou}
P(\alpha', \omega, n)= P(\alpha, \omega, n), \quad 
\Imm Z_{\omega}(\alpha ') \le \Imm Z_{\omega}(\alpha),
\end{align}
any $\omega$-Gieseker semistable 
sheaf $E$ of numerical type $\alpha '$ is $\sigma _k$-semistable.
\end{lem}
\begin{proof}
First using Lemma~\ref{useful} (i), (ii), 
the set of $\alpha' \in C(X)$ 
which satisfies (\ref{hirou}) is finite
for a fixed $\alpha$. Thus we may assume $\alpha '=\alpha$.  
Note that the 
case of $r=l=0$ is obvious. 
The case of $r>0$, $\omega \cdot l>0$ is proved 
in~\cite[Proposition 14.2]{Brs2}.
One can also check that in the proof of
 \textit{loc.cite.}, the desired $N>0$ is 
taken to be uniformly for any $\omega$-Gieseker semistable sheaf $E$
of numerical type $\alpha$.
(We leave the readers to check the detail. It is enough to 
notice 
in~\cite[Lemma 14.3]{Brs2}
that the set of $\omega$-Gieseker semistable sheaves 
of numerical type $\alpha$ is bounded.)
Thus it is enough to check the case of $r=0$ and $l\neq 0$. 
Let $E$ be a $\omega$-Gieseker semistable sheaf with 
$v(E)=(0,l,s)$. 
Since $\phi _k$ goes to $1/2$ for $k\to \infty$, we may assume 
$1/4<\phi _k < 3/4$. 
For each $k$, let $E_{k}\in \aA _{\omega}$ be the $\sigma _k$-semistable
factor of $E$ whose phase is the largest.
 If $E$ is not semistable in $\sigma _k$, 
we have 
\begin{align}\label{suzuki}
\phi _k(E_k) > \phi _k > \frac{1}{4}.
\end{align}
We have the exact sequence in $\aA _{\omega}$, 
\begin{align}\label{lala}
 0 \lr E _k \lr E \lr E_k ' \lr 0
\end{align}
Then the associated long exact sequence of (\ref{suzuki}) 
with respect to the standard t-structure 
implies that $E_k$ is a sheaf. We have the sequence, 
\begin{align}
\label{suzuki2}
 0 \lr (E_k)_{\rm{tor}} \lr E_k \lr (E_k)_{\rm{fr}} \lr 0,
\end{align}
which is exact in both
$\aA _{\omega}$ and $\Coh (X)$. 
Combining sequences (\ref{lala}) and (\ref{suzuki2}), we obtain the 
exact
sequence in $\aA _{\omega}$, 
\begin{align}\label{suzuki3}0 \lr (E_k)_{\rm{tor}} \lr E \lr F \lr 0,
\end{align}
Again the long exact sequence associated to (\ref{suzuki3}) 
implies that (\ref{suzuki3}) is also exact in $\Coh (X)$. 
Because $E$ is $\omega$-Gieseker semistable, we have
$P((E_k)_{\rm{tor}}, \omega, n)\le P(E, \omega, n)$.
Thus we have 
\begin{align}\label{suzuki4}
\phi _k((E_k)_{\rm{tor}}) \le \phi _k \le 3/4,\end{align}
by Lemma~\ref{useful} (iii). 
Then the sequence (\ref{suzuki2}) and (\ref{suzuki})
imply the map 
$$k\longmapsto \frac{\Ree Z_{k\omega}((E_k)_{\rm{fr}})}{\Imm Z_{k\omega}((E_k)_{\rm{fr}})},$$
is bounded above. Then applying Lemma~\ref{con1} and Lemma~\ref{fina}, 
there is $N>0$ such that
the map $k\mapsto \Ree Z_{\omega}((E_k)_{\rm{fr}})$ is bounded above
for $k\ge N$. 
Hence by Lemma~\ref{asu}, the set 
$$\{ v((E_k)_{\rm{fr}}) \in \NS ^{\ast}(X) \mid k\in \mathbb{Q}_{\ge N}\},$$ 
is a finite set. 
Thus we have 
 \begin{align}\label{contra}
 \phi _k ((E_k)_{\rm{fr}}) \lr 0, \end{align} for 
$k\to \infty$. 
However since we have (\ref{suzuki2}) and (\ref{suzuki4}), 
(\ref{contra}) implies that
$\phi _k(E _k) < 1/4$ for $k\ge N$ by replacing $N$ if necessary.
This contradicts to $(\ref{suzuki})$, 
 thus for such $N$ and $k\ge N$, $E$ must be
  $\sigma _k$-semistable. 
The above proof also shows that one can take $N$ 
uniformly for all $\omega$-Gieseker semistable sheaf 
$E$ of numerical type $\alpha$.

\end{proof}

Finally we show the following. 
\begin{thm}\label{goal}
For $\alpha \in C(X)$, we have $J^{\alpha}=\hat{J}^{\alpha}$. 
\end{thm}
\begin{proof}
Since $v(\alpha \otimes \lL)=v(\alpha)\cdot \ch (\lL)$ for 
$\lL \in \Pic (X)$, 
Lemma~\ref{line} implies that
we may assume $v(\alpha)=(r,l,s)$ with $\omega \cdot l>0$ or $r=l=0$. 
It is enough to compare $J^{\alpha}(\sigma _k)$ and 
$\hat{J}^{\alpha}(\omega)$ for
 $k\ge N$, where $N$ is chosen as in Proposition~\ref{semi} and 
 Lemma~\ref{koware}.
Take $\alpha _1, \cdots, \alpha _n \in C^{\sigma _k}(\phi _k)$ 
such that  $\alpha _1 + \cdots + \alpha _n =\alpha$
and $\prod _{i=1}^n I^{\alpha _i}(\sigma _k)\neq 0$. 
Then first applying Proposition~\ref{semi}, we have
$$
\alpha _i \in C(X), \quad 
\mM ^{(\alpha _i, \phi _k)}(\sigma _k) \subset \hat{\mM}^{\alpha _i}(\omega).
$$
For a fixed $k\ge N$, 
let $\sigma _k \in \mathfrak{B}^{\circ}$
 be an open neighborhood of $\sigma _k$ 
such that its closure $\mathfrak{B}$ is compact.
Then there is a wall and chamber structure 
$\{ \wW _{\gamma}\}_{\gamma \in \Gamma}$
on 
$\mathfrak{B}$ with respect to (\ref{mass}). 
There is a subset $\Gamma '\subset \Gamma$ and a connected 
component $\cC$ as in (\ref{component})
such that infinitely many $\sigma _{k'}$ for $k'\ge \mathbb{Q}_{\ge N}$ 
are contained in $\cC$. We may assume $\sigma _k \in \cC$. 
Then if $\alpha _i$ and $\alpha _j$ are not proportional 
in $\nN (X)$, we have
$$\Imm \frac{Z_{k'\omega}(\alpha _j)}{Z_{k'\omega}(\alpha _i)}=0,$$
for infinitely many $k' \in \mathbb{Q}_{\ge N}$.
By Lemma~\ref{useful} (iv), this implies 
$$P(\alpha _i, \omega,n)=P(\alpha_j, \omega,n)=P(\alpha, \omega, n),$$
for any $i,j$. Then one can apply Lemma~\ref{koware} 
and conclude
\begin{align}\label{hell}
\mM ^{(\alpha _i, \phi _k)}(\sigma _k) = \hat{\mM}^{\alpha _i}(\omega).
\end{align}
Hence we have  $\prod _{i=1}^{n} I^{\alpha _i}(\sigma _k)=
\prod _{i=1}^n \hat{I}^{\alpha _i}(\omega)$. 

Conversely take $\alpha _1, \cdots, \alpha _n \in C(X)$ such that
$\prod _{i=1}^n \hat{I}^{\alpha _i}(\omega)\neq 0$ and  
$\alpha _1 + \cdots + \alpha _n =\alpha$, 
$P(\alpha _i, \omega,n)=P(\alpha, \omega,n)$. 
Again (\ref{hell}) holds for $k\ge N$ by Proposition~\ref{semi} and 
Lemma~\ref{koware}, 
 so $\prod _{i=1}^{n} I^{\alpha _i}(\sigma _k)=
\prod _{i=1}^n \hat{I}^{\alpha _i}(\omega)$ holds. Also
 $P(\alpha _i, \omega, n)=P(\alpha, \omega, n)$ implies 
 $\alpha _i \in C^{\sigma _k}(\phi _k)$. 
 Thus the sum (\ref{yomi}) and (\ref{countsheaf}) are equal. 

\end{proof}

\begin{rmk}\emph{
In this paper, we do not give the explicit computation of the 
invariant $J^{\alpha}(\sigma)$. However if $\alpha \in C(X)$ is primitive, 
then by the work of Yoshioka~\cite{Yoshi},
$J^{\alpha}(\sigma)$ can be
 computed by the invariant of the Hilbert scheme of 
points on $X$. As commented in~\cite{Joy4}, it might be possible 
to compute the invariant for other $\alpha \in \nN (X)$
using this remark and Theorem~\ref{mainth}.} 
\end{rmk}

\bibliographystyle{jplain}
\bibliography{math}

\begin{thebibliography}{10}

\bibitem{AP}
D.~Abramovich and A.~Polishchuk.
\newblock Sheaves of t-structures and valuative criteria for stable complexes.
\newblock {\em J.reine.angew.Math}, Vol. 590, pp. 89--130, 2006.

\bibitem{Ber}
A.~Bergman.
\newblock Stability conditions and {B}ranes at {S}ingularities.
\newblock {\em preprint}.
\newblock math.AG/0702092.

\bibitem{Borch}
R.~Borcherds.
\newblock Automorphic forms on ${O}_{s+2,s}(\mathbb{R})$ and infinite products.
\newblock {\em Invent.math}, Vol. 120, pp. 161--213, 1995.

\bibitem{Brs6}
T.~Bridgeland.
\newblock Spaces of stability conditions.
\newblock {\em preprint}.
\newblock math.AG/0611510.

\bibitem{Brs3}
T.~Bridgeland.
\newblock Stability conditions and {K}leinian singularities.
\newblock {\em preprint}.
\newblock math.AG/0508257.

\bibitem{Brs2}
T.~Bridgeland.
\newblock Stability conditions on ${K}$3 surfaces.
\newblock {\em preprint}.
\newblock math.AG/0307164.

\bibitem{Brs1}
T.~Bridgeland.
\newblock Stability conditions on triangulated categories.
\newblock {\em Ann of Math (to appear)}.
\newblock math.AG/0212237.

\bibitem{ICM}
T.~Bridgeland.
\newblock Derived categories of coherent sheaves.
\newblock {\em Proceedings of the 2006 ICM}, 2006.
\newblock math.AG/0602129.

\bibitem{Brs4}
T.~Bridgeland.
\newblock Stability conditions on a non-compact {C}alabi-{Y}au threefold.
\newblock {\em Comm. Math. Phys}, Vol. 266, pp. 715--733, 2006.

\bibitem{Dou1}
M.~Douglas.
\newblock D-branes, categories and ${N}=1$ supersymmetry.
\newblock {\em J.Math.Phys}, Vol.~42, pp. 2818--2843, 2001.

\bibitem{Dou2}
M.~Douglas.
\newblock Dirichlet branes, homological mirror symmetry, and stability.
\newblock {\em Proceedings of the 1998 ICM}, pp. 395--408, 2002.
\newblock math.AG/0207021.

\bibitem{Hu}
D.~Huybrechts and M.Lehn.
\newblock {\em Geometry of moduli spaces of sheaves}, Vol. E31 of {\em Aspects
  in Mathematics}.
\newblock Vieweg, 1997.

\bibitem{Inmo}
M.~Inaba.
\newblock Moduli of stable objects in a triangulated category.
\newblock {\em preprint}.
\newblock math.AG/0612078.

\bibitem{Inaba}
M.~Inaba.
\newblock Toward a definition of moduli of complexes of coherent sheaves on a
  projective scheme.
\newblock {\em J.Math.Kyoto Univ.}, Vol. 42-2, pp. 317--329, 2002.

\bibitem{IUU}
A.~Ishii, K.Ueda, and H.Uehara.
\newblock Stability {C}onditions on ${A}_n$-{S}ingularities.
\newblock math.AG/0609551.

\bibitem{Joy3}
D.~Joyce.
\newblock Configurations in abelian categories
  {I}\hspace{-.1em}{I}\hspace{-.1em}{I}. {S}tability conditions and identities.
\newblock {\em preprint}.
\newblock math.AG/0410267.

\bibitem{Joy4}
D.~Joyce.
\newblock Configurations in abelian categories {I}\hspace{-.1em}{V}.
  {I}nvariants and changing stability conditions.
\newblock {\em preprint}.
\newblock math.AG/0410268.

\bibitem{Joy5}
D.~Joyce.
\newblock Motivic invariants of {A}rtin stacks and `stack functions'.
\newblock {\em preprint}.
\newblock math.AG/0509722.

\bibitem{Joy1}
D.~Joyce.
\newblock Configurations in abelian categories {I}. {B}asic properties and
  moduli stack.
\newblock {\em Advances in Math}, Vol. 203, pp. 194--255, 2006.

\bibitem{Joy2}
D.~Joyce.
\newblock Configurations in abelian categories {I}\hspace{-.1em}{I}.
  {R}ingel-{H}all algebras.
\newblock {\em Advances in Math}, Vol. 210, pp. 635--706, 2007.

\bibitem{Joy}
D.~Joyce.
\newblock Holomorphic generating functions for invariants counting coherent
  sheaves on {C}alabi-{Y}au 3-folds.
\newblock {\em Geometry and Topology}, Vol.~11, pp. 667--725, 2007.

\bibitem{GL}
G.~Laumon and L.~Moret-Bailly.
\newblock {\em Champs alg{\'e}briques}, Vol.~39 of {\em Ergebnisse der
  Mathematik und ihrer Grenzgebiete}.
\newblock Springer Verlag, Berlin, 2000.

\bibitem{LIE}
M.~Lieblich.
\newblock Moduli of complexes on a proper morphism.
\newblock {\em J.Algebraic Geom}, Vol.~15, pp. 175--206, 2006.

\bibitem{Mac}
E.~Macri.
\newblock Some examples of moduli spaces of stability conditions on derived
  categories.
\newblock {\em preprint}.
\newblock math.AG/0411613.

\bibitem{MW}
K.~Matsuki and R.Wentworth.
\newblock Mumford-{T}haddeus principle on the moduli space of vector bundles on
  an algebraic surface.
\newblock {\em Internat.J.Math}, Vol.~8, pp. 97--148, 1997.

\bibitem{Oka}
S.~Okada.
\newblock Stability manifold of $\mathbb{P}^1$.
\newblock {\em J.Algebraic Geom}, Vol.~15, pp. 487--505, 2006.

\bibitem{Tho}
R.~Thomas.
\newblock Stability conditions and the braid groups.
\newblock {\em Comm.Anal.Geom}, Vol.~14, pp. 135--161, 2006.

\bibitem{Tst2}
Y.~Toda.
\newblock Stability conditions and {C}alabi-{Y}au fibrations.
\newblock {\em preprint}.
\newblock math.AG/0068495.

\bibitem{Tst}
Y.~Toda.
\newblock Stability conditions and crepant small resolutions.
\newblock {\em Trans.Amer.Math.Soc (to appear)}.
\newblock math.AG/0512648.

\bibitem{Yoshi}
K.~Yoshioka.
\newblock Moduli spaces of stable sheaves on abelian surfaces.
\newblock {\em Math.Ann}, Vol. 321, pp. 817--884, 2001.

\end{thebibliography}

Yukinobu Toda, Graduate School of Mathematical Sciences, University of Tokyo

\textit{E-mail address}:toda@ms.u-tokyo.ac.jp

\end{document}